\documentclass[12pt]{amsart}

\extrafloats{100}

\newcommand{\br}{\textrm{br}}
\newcommand{\Sc}{\mathcal{S}}
\newcommand{\Tc}{\mathcal{T}}
\newcommand{\Hc}{\mathcal{H}}
\newcommand{\C}{\mathbb{C}}
\newcommand{\CP}{\mathbb{P}_\mathbb{C}}
\newcommand{\CH}{\mathbb{H}_\mathbb{C}}
\newcommand{\B}{\mathbb{B}}
\newcommand{\R}{\mathbb{R}}
\newcommand{\Q}{\mathbb{Q}}
\newcommand{\F}{\mathbb{F}}
\newcommand{\N}{\mathbb{N}}
\newcommand{\Z}{\mathbb{Z}}
\newcommand{\IK}{\mathbb{K}}

\newcommand{\Oc}{\mathcal{O}}

\usepackage{amssymb}
\usepackage{a4wide}
\usepackage{url}
\usepackage{alltt}

\newcommand{\tr}{\textrm{tr}}
\newcommand{\inn}[2]{\langle #1,#2\rangle}
\newcommand{\FPG}[2]{\left\langle\ #1\ |\ #2\ \right\rangle} 

\newtheorem{thm}{Theorem}[section]
\newtheorem{lem}{Lemma}[section]
\newtheorem{dfn}{Definition}[section]
\newtheorem{prop}[thm]{Proposition}
\newtheorem{cor}[thm]{Corollary}
\theoremstyle{remark}
\newtheorem{rmk}[thm]{Remark}

\def\cqfd{\mbox{}\nolinebreak\hfill$\Box$\medbreak\par}
\newenvironment{pf}{\noindent\textbf{Proof:}}{\cqfd}

\title{On subgroups of finite index in complex hyperbolic lattice triangle groups}

\author{Martin Deraux}

\address{Martin Deraux : Universit\'e Grenoble Alpes, Institut
  Fourier, 100 rue des Math\'ematiques 38610 Gi\`eres; Sorbonne
  Université and Université de Paris, CNRS, INRIA, IMJ-PRG, Ouragan,
  F-75005 Paris, France } \email{martin.deraux@univ-grenoble-alpes.fr}

\date{June 20, 2023}

\begin{document}

\begin{abstract}
  We study several explicit finite index subgroups in the known
  complex hyperbolic lattice triangle groups, and show some of them
  are neat, some of them have positive first Betti number, some of
  them have a homomorphisms onto a non-Abelian free group.  For some
  lattice triangle groups, we determine the minimal index of a neat
  subgroup. Finally, we answer a question raised by Stover and
  describe an infinite tower of neat ball quotients all with a single
  cusp.
\end{abstract}

\maketitle

\section{Introduction}

The goal of this paper is to study some explicit finite index
subgroups of the known complex hyperbolic lattice triangle
groups. Recall that the first examples of such lattices were studied
by Mostow~\cite{mostow-pacific} (see
also~\cite{deligne-mostow},~\cite{thurston-shapes} and the references
given there). New examples were discovered quite a bit later in joint
work of the author with Parker aud Paupert,
see~\cite{dpp1},~\cite{dpp2}.

Even though triangle groups are of course very special among all
lattices in $PU(2,1)$, they allow us to describe all commensurability
classes of non-arithmetic lattices in $PU(2,1)$ that are known to this
day, so we consider it an important and interesting class.

In the papers cited above, explicit fundamental domains for the action
of complex hyperbolic triangle groups are described, which allows us
to describe
\begin{itemize}
\item a presentation for these groups in terms of generators and
  relations;
\item the conjugacy classes of isotropy groups;
\item the conjugacy classes of cusps.
\end{itemize}
Even though most of that information can be gathered in~\cite{dpp2},
that task would require using the output of our computer program
Spocheck~\cite{spocheck}, which is probably too much to ask from the
average reader.

In this paper we explain some details on how to achieve the last two
items, we give enough information for the above data to be
recontructed fairly convieniently, and we list the results in the form
of a computer file~\cite{computer-code}. We also give some
applications to the construction of explicit finite index subgroups of
(some of) the complex hyperbolic lattice triangle groups.

We focus on three different kinds of subgroups, namely
\begin{itemize}
\item torsion-free/neat subgroups;
\item subgroups with positive first Betti number;
\item subgroups that surject onto a non-abelian free group.
\end{itemize}

Recall that a lattice $\Gamma\subset PU(2,1)$ is torsion-free
(i.e. contains is no non-trivial element of finite order) if and only
if the quotient map $\CH^2\rightarrow \CH^2/\Gamma$ is an unramified
covering. Equivalently, one requires that the quotient $\CH^2/\Gamma$
is a manifold with charts given by local inverses of the quotient map
from $\CH^2$.

A lattice is neat (see~\cite{borel}) if it is torsion-free, and every
parabolic element in the group can be realized by a unipotent
matrix. Neat lattices are also important, because when a lattice
$\Gamma\subset PU(2,1)$ is neat, the quotient $X=\CH^2/\Gamma$ admits a
smooth toroidal compactification, obtained by adding to each end of
$X$ an elliptic curve with negative self-intersection (see~\cite{amrt}
or~\cite{mok-projective}).

A well-known result of Selberg, usually referred to as Selberg's lemma
(see~\cite{selberg},~\cite{alperin}), says that every finitely
generated matrix group admits a neat subgroup of finite
index. Selberg's argument relies on taking a principal congruence
subgroup modulo a well-chosen ideal, and this tends to produce
subgroups of very large index (see the tables of
section~\ref{sec:neat} for examples). We will attempt to get
torsion-free subgroups of reasonably small index, so that one can hope
to get a presentation, compute the abelianization, count and describe
the cusps.

In the following discussion, let $\Gamma$ stand for a lattice triangle
group. By construction $\Gamma$ is not torsion-free, it is in fact
generated by a special kind of torsion element (namely complex
reflections, which are elliptic transformations with a repeated
eigenvalue, see section~\ref{sec:chg}).

It is an standard consequence of the orbit-stabilizer theorem that the
index of a torsion-free subgroup $H\subset \Gamma$ must be a multiple
of the least common multiple $L$ of the orders of isotropy groups in
$\Gamma$ (see~\cite{eek} for example). Moreover, if $\Gamma$ (hence
also $H$) is cocompact, a standard consequence of Noether's formula
(see Proposition~\ref{prop:euler-three}) implies that
$\chi(H\backslash\CH^2)$ must be a multiple of 3, which in some cases
gives a slightly larger lower bound for the index $[\Gamma:H]$.

We call the \emph{obvious lower bound} the number
\begin{equation}
  \label{eq:opt}
  L^{opt}=\left\{\begin{array}{l}
                   3L \textrm{ if } \Gamma \textrm{ is cocompact and } 3\chi^{orb}(\CH^2/\Gamma)\notin 3\Z\\
                   L \textrm{ otherwise}\end{array}\right.
\end{equation}
If there exists a neat subgroup of index
$L^{opt}$, we will say that {\bf the obvious lower bound is realized}.

For lattices in $PU(1,1)\cong PSL(2,\R)$, the analogue of the obvious
lower bound (with a factor of 2 instead of 3) is realized for every
lattice, by a result of Edmonds, Ewing and Kulkarni~\cite{eek}. For
$PSL(2,\C)$, Jones and Reid~\cite{jones-reid} have shown that there
are lattices where the obvious lower bound is arbitrarily far from
being realized.

About lattice triangle groups in
$PU(2,1)$, we prove the following using the results in this paper in
conjunction with Magma calculations, see section~\ref{sec:neat}.
\begin{thm}\label{thm:main}
  The obvious lower bound is realized for 
  $\Sc(4,\sigma_1)$, 
  $\Sc(3,\sigma_5)$,
  $\Sc(5,\sigma_{10})$,
  $\Sc(10,\sigma_{10})$,
  $\Tc(3,{\bf S_2})$,
  $\Tc(3,{\bf E_2})$,
  $\Tc(4,{\bf E_2})$,
  $\Tc(5,{\bf H_2})$,
  $\Gamma(4,5/12)$,
  $\Gamma(6,1/3)$,
  $\Gamma(3,1/3)$,
  $\Gamma(4,1/4)$,
  $\Gamma(5,1/10)$,
  $\Gamma(10,0)$,
  $\Gamma(3,1/6)$,
  $\Gamma(4,1/12)$.
\end{thm}
\begin{cor}
  The minimal index of torsion-free subgroups in the lattice
  triangle groups in the statement of Theorem~\ref{thm:main} are as in Table~\ref{tab:opt-index}.\\
  \begin{table}
    \centering
    \begin{tabular}[htbp]{|c|c|c|c|c|c|c|c|}
    \hline
    $\Sc(4,\sigma_1)$ & $\Sc(3,\sigma_5)$ & $\Sc(5,\sigma_{10})$ & $\Sc(10,\sigma_{10})$ & $\Tc(3,{\bf S_2})$ & $\Tc(3,{\bf E_2})$ & $\Tc(4,{\bf E_2})$ & $\Tc(5,{\bf H_2})$\\
    \hline
    $96$              &   $360$           &     $600$            &    $300$              &    $360$           &  $72$              &   $96$             &   $600$\\
    \hline
  \end{tabular}\\
  \begin{tabular}[htbp]{|c|c|c|c|c|c|c|c|c|}
    \hline
    $\Gamma(4,\frac{5}{12})$ & $\Gamma(6,\frac{1}{3})$ & $\Gamma(3,\frac{1}{3})$ & $\Gamma(4,\frac{1}{4})$ & $\Gamma(5,\frac{1}{10})$ & $\Gamma(10,0)$ & $\Gamma(3,\frac{1}{6})$ & $\Gamma(4,\frac{1}{12})$ \\
    \hline
    $864$                    &   $18$                  &     $864$               &   $96$                  &    $600$                 &  $150$         &   $72$                  &   $864$\\
    \hline
  \end{tabular}
  \caption{Minimal index of a torsion-free subgroup}  \label{tab:opt-index}
  \end{table}
\end{cor}

For the other lattice triangle groups, we were not able to prove that
the obvious lower bound cannot be realized (but we suspect that this
is the case). Note that for some lattices, we get a neat finite index
subgroup of index reasonably close to $L^{opt}$; for others, the only
neat subgroups we know are principal congruence subgroups, and in some
cases these congruence subgroups seem inaccessible to computation
using current computer technology (see Tables~\ref{tab:neat-3s1}
through~\ref{tab:neat-m3-0} for the list of subgroups we were able to
obtain).

We point out that some of the results in Theorem~\ref{thm:main} are
not new.  For the group $\Gamma(3,1/3)$, which is also the group
Deligne-Mostow group $\Gamma_{\mu,S_3}$ for $\mu=(2,2,2,7,11)/12$, the
subgroup of index 864 with Abelianization $\Z^2$ is the fundamental
group of the Cartwright-Steger surface.

For the groups $\Gamma(6,1/3)$ and $\Gamma(3,1/6)$, the subgroups of
optimal index give a neat subgroup $H$ with
$\chi(H\backslash\CH^2)=1$, 4 cusps and Abelianization $\Z^4$, which
correspond to Hirzebruch's surface
(see~\cite{hirzebruch},~\cite{dicerbo-stover}).

Note also that there are 4 non-arithmetic lattices in the list given
in Theorem~\ref{thm:main}, namely $\Sc(3,\sigma_1)$,
$\Sc(3,\sigma_5)$, $\Tc(4,{\bf E_2})$, $\Gamma(4,1/12)$. For these
lattices, getting a subgroup with optimal index yields a ball quotient
with fairly large Euler characteristic (namely $\chi=42,98,102,39$
respectively).

We point out a nice by-product of our search for neat subgroups in
lattice triangle groups, which answers a question raised by Stover
in~\cite{stover-cusps} (the analogous question for real hyperbolic
4-manifolds was raised by Long-Reid~\cite{long-reid}, and
answered by Kolpakov and Martelli~\cite{kolpakov-martelli}).
\begin{thm}\label{thm:one-cusp}
  The Mostow group $G=\Gamma(6,0)$ has a neat subgroup $H$ of index 72
  such that $X=H \backslash \CH^2$ has exactly one cusp. 
\end{thm}
As pointed out to us by Matthew Stover, this result can be improved to
the following statement.
\begin{thm}\label{thm:tower}
  There exists an infinite tower
  $H=H_1\supset H_2 \supset H_3\supset \cdots$ of neat subgroups with
  $[H:H_n]\rightarrow\infty$ such that for every $n$,
  $X_n=H_n\backslash\CH^2$ has exactly one cusp.
\end{thm}
A similar infinite tower of neat subgroup with 2 cusps was
constructed by Di Cerbo and Stover~\cite{dicerbo-stover-if}.

About subgroups with positive first Betti numbers (see
section~\ref{sec:betti}), first note that we record the first Betti
number of the torsion-free subgroups in the last column of
Tables~\ref{tab:neat-3s1} through~\ref{tab:neat-m3-0}. A lot of these
have $b_1>0$ (this is equivalent to saying that their Abelianization
is infinite). In cases where all the torsion-free subgroups we have
found have $b_1=0$, we try computing all (normal) subgroups of
``small'' index, and check if these subgroups have $b_1>0$; the word
``small'' has no precise meaning here, we simply mean that the Magma
command for computing the corresponding subgroups should take less
than a day, say.  For several lattice triangle groups, both of these
methods fail, and we do not know any explicit subgroup of finite index
with $b_1>0$.

Finally we review some facts (that are probably well-known to experts)
that show that many Mostow groups are large, i.e. they admit a finite
index subgroup that maps onto a non-Abelian free group. The subgroups
in question are obtained by constructing explicit subgroups that
admit a homomorphism onto a hyperbolic triangle group. Largeness then
follows from largeness of every Fuchsian group (see
Lemma~\ref{lem:sl2r}).

\noindent
\textbf{Acknowledgements:} I wish to thank Philippe Eyssidieux and
Julien Paupert for interesting discussions related to this work. Warm
thanks also go to Matthew Stover, for patiently explaining a lot of
the computational techniques used in the paper, and for suggesting to
strengthen Theorem~\ref{thm:one-cusp} to~\ref{thm:tower}. Finally, the
author acknowledges support from INRIA, in the form of a research
semester in the ``Ouragan'' team.

\section{Background on lattice triangle groups}

\subsection{The complex hyperbolic plane}\label{sec:chg}

We briefly review basic facts about the complex hyperbolic plane and
complex hyperbolic triangle groups.  Recall that the complex
hyperbolic plane $\CH^2$ is the only Hermitian symmetric space of
complex dimension 2 with constant negative holomorphic sectional
curvature (see~\cite{kobayashi-nomizu}). One possible model is the
unit ball in $\B\subset \C^2$ equipped with the Bergman metric, but it
is convenient to see $\C^2$ as an affine chart of $\CP^2$, and to work
in homogeneous coordinates. As a set, $\CH^2$ is the set of complex
lines that are negative for a fixed Hermitian inner product
$\inn{\cdot}{\cdot}$ of signature $(2,1)$ on $\C^{3}$, which we can
describe as $\inn{V}{W}=W^*HV$ for some Hermitian matrix $H$. Recall
that $U(H)$ denotes the group
$$
U(H)=\{A\in GL_{3}(\C):A^*HA=H\},
$$
and $PU(H)$ is the quotient of $U(H)$ by the subgroup of scalar
matrices. Note that by Sylvester's law of inertia, up to isomorphism,
the real Lie group $PU(H)$ is independent of the choice of Hermitian
form $H$, a common choice being the diagonal matrix
$\textrm{diag}(-1,1,1)$.

Up to scale, there is a unique $PU(H)$-invariant Riemannian metric on
$\CH^2$, which is K\"ahler and has constant holomorphic sectional
curvature. We refer to $\CH^2$, equipped with the unique
$PU(H)$-invariant Riemannian metric of holomorphic sectional curvature
$-1$, as the complex hyperbolic plane. When we talk about volume, we
mean volume with respect to the corresponding Riemannian volume
form.

The stabilizer of a point $x=[V]$ in $\CH^2$ in $PU(H)$ is isomorphic
to group of isometries of the restriction of the Hermitian form to
$V^\perp=\{W\in\C^3:\inn{V}{W}=0\}$. Because of the fact that $V$ is a
negative vector and $H$ has signature $(2,1)$, $V^\perp$ is a positive
definite subspace, so the stabililizer of a point is isomorphic to
$U(2)$.

One checks that $PU(H)$ is the full group of homomorphic isometries,
and that it has index two in the full group of isometries (an isometry
not in $PU(H)$ is induced by $V\mapsto Q\bar V$ where $Q$ is a matrix
such that $Q^*\bar H Q=H$ (such a matrix exists because $H$ and
$\bar H$ have the same signature).

An explicit expression for the invariant Riemannian metric can be
found in~\cite{kobayashi-nomizu}, see also~\cite{goldman-book}. We
will not need any expression for that metric; for concreteness, we
describe the corresponding distance function, which is given by
$$
\cosh\left(\frac{1}{2}d(\C V,\C W)\right)=\frac{|\inn{V}{W}|}{\sqrt{\inn{V}{V}\inn{W}{W}}}.
$$
where $V,W\in\C^3$ are any two negative vectors. The factor
$\frac{1}{2}$ is included for the constant holomorphic sectional
curvature to be equal to $-1$.

Given two distinct points $\C V,\C W\in\CH^2$, i.e. $V,W$ are linearly
independent over $\C$, the restriction of the Hermitian form to the
complex span $\C\{V,W\}$ has signature $(1,1)$, so the set of negative
vectors in the span gives a copy of $\CH^1$, which is totally
geodesic. We call this {\bf the complex line through $V$ and $W$}.

\subsection{Point stabilizers}

Recall that the isometric action of $PU(H)$ on the open ball $\CH^2$
extends to an action on the closed ball
$\CH^2\cup\partial_\infty\CH^2$ by homeomorphisms. Every homeomorphism
of a closed ball has a fixed point, and the position of the
corresponding fixed points (in $\CH^2$ or in $\partial_\infty\CH^2$)
gives a classification of isometries into elliptic, parabolic,
loxodromic isometries (see section~6.2 in~\cite{goldman-book}).

Elliptic isometries are the ones that have a fixed point in
$\CH^2$. The ones without any fixed point in $\CH^2$ are either
parabolic (unique fixed point in $\partial_\infty\CH^2$) or loxodromic
(two distinct fixed points in $\partial_\infty\CH^2$).

Suppose $A\in PU(H)$ is elliptic. Then $A$ has an eigenvector $V_0$
with $\inn{V_0}{V_0}<0$, and by multiplying $A$ by a constant, we may
assume the corresponding eigenvalue is $1$. Moreover, we may normalize
$V_0$ so that $\inn{V_0}{V_0}=-1$.  The orthogonal complement
$V_0^\perp=\{W\in\C^3:\inn{V_0}{W}=0\}$ is invariant under $A$ and,
because of the signature assumption, the restriction of the Hermitian
form to $V^\perp$ is positive definite. Pick a $V_1\in V_0^\perp$,
which once again we may normalize so that $\inn{V_1}{V_1}=1$, and then
take $V_2$ to be the Hermitian box product $V_0\boxtimes V_1$
(i.e. the usual cross-product of $V_0^*H$ with $V_1^*H$). Normalizing
$V_2$, we get $A$ to be diagonal in the basis $V_0,V_1,V_2$, which is
standard Lorentzian in the sense that the matrix of the Hermitian form
in that basis is
$$
J=\left(\begin{matrix}
 -1 & 0 & 0\\
  0 & 1 & 0\\
  0 & 0 & 1
  \end{matrix}\right),
$$
As mentioned before, this identifies the stabilizer of $\C V_0$ in
$PU(H)$ with the unitary group $U(2)$ (which is a compact group).

When $A$ has distinct eigenvalues, $A$ is called a {\bf regular
  elliptic} isometry. In that case, $A$ has three linearly independent
eigenvectors, and the Lorentzian character of the diagonalizing basis
is essentially automatic because eigenvectors with different
eigenvalues are orthogonal with respect to $\inn{\cdot}{\cdot}$ (one
simply normalizes any triple of eigenvectors).

When $A$ is elliptic and has a double eigenvalue, the above argument
shows that $A$ remains diagonalizable (because every matrix in $U(2)$
is diagonalizable). We can write
$$
\C^3=\C V_0\oplus V_0^{\perp}
$$
for some $V_0\in\C^3$ with $\inn{V_0}{V_0}\neq 0$ (if the inner
product were $0$, the above decomposition would not be a direct
sum). Note that $A$ acts projectively as the identity on the
projectivization of $V_0^\perp$, because it is an eigenspace for $A$.

If $\inn{V_0}{V_0}>0$, we call $A$ a {\bf complex reflection}, and we
call (the projectivization of) $V_0^\perp$ its {\bf mirror}, which is
a complex line, i.e.  a totally geodesic copy of $\CH^1$ in $\CH^2$.

If $\inn{V_0}{V_0}<0$, then $\C V_0$ gives an isolated fixed point in
$\CH^2$ for the action of $A$, and we say that $A$ is a {\bf complex
reflection in a point}.

In both cases there is a simply formula for $A$ in terms of the
Hermitian inner product. Indeed, let $\zeta\in\C$ be such that
$|\zeta|=1$ and consider the linear map given for every $X\in\C^3$ by
\begin{equation}\label{eq:refl}
R(X)=X+(\zeta-1)\frac{\inn{X}{V_0}}{\inn{V_0}{V_0}}V_0.
\end{equation}
It is easy to see that $V_0^\perp$ and $\C V_0$ are eigenspaces of
$R$, with eigenvalue $1$ and $\zeta$ respectively. Depending on the
sign of $\inn{V_0}{V_0}$, we get either a complex reflection (with
mirror $V_0^\perp$) or a complex reflection in a point. The complex
number $\zeta$ is called the {\bf multiplier} of the complex
reflection.

\begin{rmk}
  \begin{enumerate}
  \item It follows from the above discussion that if a holomorphic
    isometry $A\in PU(H)$ fixes two distinct points in $\CH^2$, then
    it fixes the complex lines through these two points, hence it is a
    complex reflection.
  \item Regular elliptic elements can have non-trivial powers that are
    no longer regular.
  \end{enumerate}
\end{rmk}

\subsection{Ideal point stabilizers}\label{sec:heisenberg}

We now review some facts about the stabilizer in $G=PU(H)$ of an ideal
point, i.e. a complex line spanned by a null vector. Recall that all
Hermitian forms of signature $(2,1)$ on $\C^3$ are equivalent over
$\R$, so we may assume by choosing a suitable basis of $\C^3$ that the
matrix of the Hermitian form in the standard basis $e_1,e_2,e_3$ of
$\C^3$ is the antidiagonal matrix
$$
J=\left(\begin{matrix}
  0 & 0 & 1\\
  0 & 1 & 0\\
  1 & 0 & 0
  \end{matrix}\right).
$$
With such normalization, $e_1$ is a null vector, so
$\C e_1\in\partial_\infty\CH^2$.

The unipotent stabilizer of $e_1$ in $G$ forms a subgroup of
$Stab_G(e_1)$, which is given by the matrices of the form
\begin{equation}\label{eq:heis-tsl}
T(z,t)=\left(
  \begin{matrix}
    1 & -\bar z &\frac{-|z|^2+it}{2}\\
    0 &    1    &    z\\
    0 &    0    &    1
  \end{matrix}\right).
\end{equation}
where $z\in\C$, $t\in\R$. Let $V=(v_1,v_2,v_3)\in\C^3$ be any null
vector for the form $\inn{V}{W}=W^*JV$, i.e.
$2\Re(v_1\bar v_3)+|v_2|^2=0$. If $v_3=0$, then $v_2=0$ and $V$ is a
multiple of $e_1$; if $v_3\neq 0$, we can normalize the homogeneous
coordinates so that $v_3=1$, and then we have
$2\Re(v_1)+|v_2|^2=0$. This implies that $v$ is in the image of
$e_3=(0,0,1)$ under some transformation $T(z,t)$ (take $z=v_2$,
$t=(2v_1+|v_2|^2)/i$), i.e. the unipotent stabilizer acts transitively
on $\partial_\infty\CH^2\setminus \C e_1$.

In fact, one verifies that $T(z,t)T(z',t')=T((z,t)\star(z',t'))$,
where $\star$ denotes the Heisenberg group law, namely
$$
(z,t)\star(z',t')=(z+z',t+t'+2\Im(z\bar z')).
$$
The group $\C\times\R$ equipped with the Heisenberg group law is the
(3-dimensional real) Heisenberg group, and we denote it by $\Hc$. Note
that the center of $\Hc$ is the subgroup
$\left\{(0,t),t\in\R\right\}$.

Equation~\eqref{eq:heis-tsl} can be thought of as embedding $\Hc$ as
the unipotent stabilizer of $e_1$, or as an identification between
$\partial_\infty\CH^2\setminus\C e_1$ and $\Hc$ (the point $(z,t)$
corresponds to the complex line spanned by
$(\frac{-|z|^2+it}{2},z,1)$). Since $T(z,t)$ acts on $\Hc$ via left
$\star$-multiplication by $(z,t)$, we often refer to the
transformations $T(z,t)$ as a {\bf Heisenberg translations}. The
elements with $z=0$, which are central in the unipotent stabilizer,
are called {\bf vertical translations}.

There are complex reflections in the stabilizer, given by the diagonal
matrices $R_\zeta=\textrm{diag}(1,\zeta,1)$, $|\zeta|=1$. These act on
$\Hc$ as $(z,t)\mapsto (\zeta z,t)$, and we refer to them as {\bf
  Heisenberg rotations}. Finally, {\bf Heisenberg dilations} are given
by $\textrm{diag}(\lambda,1,1/\lambda)$ for $\lambda\in\R^*$.

The full stabilizer of $e_1$ in $G$ is generated by these three
classes (see section 4.2.2 of~\cite{goldman-book} for instance). In
terms of the rough classification of isometries into elliptic,
parabolic and loxodromic transformation, Heisenberg translations
(resp. rotations, dilations) are parabolic (res. elliptic,
loxodromic).

Note also that for arbitrary $z\in\C$, $t\in\R$ and $\zeta\in\C$ with
$|\zeta|=1$, $R_\zeta T(z,t)$ is parabolic, but it is unipotent if and
only if $\zeta=1$. The transformations $R_\zeta T(z,t)$ with
$\zeta\neq 1$ are often called {\bf screw-parabolic elements}.

\section{Lattices in $PU(2,1)$}\label{sec:lattices}

Throughout this section, $\Gamma$ denotes a lattice in $PU(H)$, in the
sense of Definition~\ref{dfn:lattice}.
\begin{dfn}\label{dfn:lattice}
  A subgroup $\Gamma\subset PU(H)$ is called a lattice if it is
  discrete and $\Gamma\backslash\CH^2$ has finite volume. A lattice is
  called cocompact (or uniform) if $\Gamma\backslash\CH^2$ is compact.
\end{dfn}
To simplify notation, we write $G=PU(H)$. Even though we will not use
this in the present paper, we point out that the lattice condition is
equivalent to the requirement that $\Gamma\backslash G$ has a finite
Haar measure (this is a consequence of the compactness of
$K=Stab_G(x)$, $x\in\CH^2$).

Since the stabilizer in $G$ of a point $x\in\CH^2$ is compact,
$Stab_\Gamma(x)$ is a finite group, for every $x\in\CH^2$. It follows
that the quotient $\Gamma\backslash\CH^2$ acquires the structure of a
complex hyperbolic orbifold.

In order for the quotient map to induce a complex hyperbolic manifold
structure on the quotient, one needs the point-stabilizers to be
trivial; since any element of finite order in $G$ fixes at least one
point in $\CH^2$, we have the following.
\begin{prop}\label{prop:TF}
  The local inverses of the quotient map
  $\CH^2\rightarrow \Gamma\backslash\CH^2$ define a manifold structure
  on $\Gamma\backslash\CH^2$ if and only if $\Gamma$ is torsion-free.
\end{prop}
When the conditions of Proposition~\ref{prop:TF} are satisfied, we say
that $\Gamma\backslash\CH^2$ is a smooth ball quotient (of finite
volume).

\begin{rmk}
  The quotient $\Gamma\backslash\CH^2$ can be a smooth even in cases
  where $\Gamma$ has non-trivial torsion; indeed, some of the ball
  quotients constructed by Hirzebruch (see~\cite{bhh},~\cite{tretkoff}
  or~\cite{deligne-mostow}) give lattices such that
  $\Gamma\backslash\CH^2$ biholomorphic to $\CP^2$ (possibly blown up
  at some points), but the corresponding quotient map
  $\CH^2\rightarrow\CP^2$ is then a branched covering.
\end{rmk}

If $\Gamma$ is torsion-free and cocompact, $X=\Gamma\backslash\CH^2$
is a compact complex surface, and Noether's formula
(see~\cite{griffiths-harris}, p.438) says that
\begin{equation}
  \Z\ni\chi(\Oc_X) = \frac{c_1(X)^2+c_2(X)}{12}.
  \label{eq:noether}
\end{equation}
Moreover, by Hirzebruch proportionality~\cite{hirzebruch-prop}, smooth
compact ball quotients satisfy the same equality of Chern numbers as
$\CP^2$, namely $c_1(X)^2=3c_2(X)$. Since $c_2(X)$ is the topological
Euler characteristic $\chi(X)$, equation~\eqref{eq:noether} implies
the following.
\begin{prop}\label{prop:euler-three}
  Let $X$ be a compact smooth ball quotient. Then $\chi(X)\in 3\Z$.
\end{prop}

There exist many compact smooth ball quotients with $\chi(X)=3$, namely
fake projective planes (these are compact ball quotients of dimension
2 with the same Betti numbers as $\CP^2$, in particular they have
$b_1=0$). Klingler~\cite{klingler} has shown that every fake
projective plane must be arithmetic (see also~\cite{yeung}), and this
sets ground for a classification, since arithmetic groups can be
described by explicit number-theoretical data
(see~\cite{weil-involutions},~\cite{tits}). Fake projective planes
were indeed enumerated by Prasad, Yeung~\cite{prasad-yeung},
Cartwright and Steger~\cite{cartwright-steger}. As a by-product of the
work needed for the classification, Cartwright and Steger also found a
compact smooth ball quotient of dimension 2 with minimal volume (this
is equivalent to having Euler characteristic 3) but which is not a
fake projective plane (in fact it has positive first Betti
number). The latter is listed in our paper, its fundamental group is
(conjugate to) a subgroup of a specific Mostow ball quotient (see
Table~\ref{tab:neat-m3-1.3}).

Note also that every element in $3\Z$ is the Euler characteristic of
some smooth compact ball quotient. Indeed, the fundamental group of
the Cartwright-Steger $\Gamma$ has Abelianization $\Z\oplus\Z$, hence
a surjective morphism $\Gamma\rightarrow \Z$; the preimage of $n\Z$ under
this morphism gives a subgroup of index $n$ in $\Gamma$, which has
Euler characteristic $3n$.

For non-compact ball quotients, the Euler characteristic can be equal
to any natural number $n\in\N^*$. This follows from the fact that the
surface constructed by Hirzebruch in~\cite{hirzebruch} has Euler
characteristic $1$ and the Abelianization of its fundamental group is
$\Z^4$, so the corresponding lattice has a homomorphism onto $\Z$ (and
the latter has subgroups of any arbitrary finite index). In fact a
stronger result was proved by Di Cerbo and
Stover~\cite{dicerbo-stover-if}, namely for every $n\in\N^*$, there
exists a non-cocompact ball quotient \emph{with two cusps} and Euler
characteristic $1$. In this paper, we will improve this to get a tower
of smooth ball quotients with a single cusp, see
Theorem~\ref{thm:tower}.

It is natural to look for a condition analogous to the torsion-free
condition but for stabilizers of ideal points. Of course a finite
covering of a non-compact ball quotient will remain non-compact, so
one cannot hope to get rid of cusps alltogether in a subgroup of
finite index (see Proposition~\ref{prop:par-lox}); the relevant
condition for cusps to be torsion-free is slightly more subtle.

We start by observing that if $\Gamma\subset G$ is discrete and
$x\in\partial_\infty\CH^2$, then $Stab_\Gamma(x)$ cannot contain both
parabolic and loxodromic elements, which simplifies the description of
stabilizers of ideal points in discrete groups, see section
~\ref{sec:heisenberg}.

Because of the lattice assumption, the thick-thin decomposition (see
section~4 of~\cite{kapovich} for instance) implies that the quotient
$\Gamma\backslash\CH^2$ has finitely many ends, each being
homeomorphic to $\mathcal{N}\times ]0,\infty[$, where $\mathcal{N}$ is
a compact quotient of the Heisenberg group, which is isomorphic to the
stabilizer in $\Gamma$ of a given ideal point
$x\in\partial_\infty\CH^2$ (see also~\cite{garland-raghunathan}). The
end of the quotient is usually referred to as a \emph{cusp}; cusps are
in 1-1 correspondence with $\Gamma$-conjugacy classes of non-trivial
isotropy groups of ideal boundary points that contain at least one
parabolic element. By extension, these stabilizers are often called
cusps as well.

In particular, we have the following.
\begin{prop}\label{prop:par-lox}
  Let $\Gamma$ be a lattice in $G=PU(H)$. The quotient
  $\Gamma\backslash\CH^2$ is cocompact if and only if $\Gamma$
  contains no parabolic element.
\end{prop}

Using the description of $Stab_G(x)$ for $x\in\partial_\infty\CH^2$
(see section~\ref{sec:chg}), we can say a little more about cusp
groups. Once again, if $Stab_\Gamma(x)$ is a cusp group, then it
consists of screw-parabolic elements, i.e. matrices of the form
\begin{equation}\label{eq:screw}
\left(
  \begin{matrix}
    1 & -\bar z &\frac{-|z|^2+it}{2}\\
    0 &  \zeta  &    z\\
    0 &    0    &    1
  \end{matrix}\right),
\end{equation}
$z,\zeta\in\C,t\in\R$, $|\zeta|=1$.

The projection
\begin{equation}
\left(
  \begin{matrix}
    1 & -\bar z &\frac{-|z|^2+it}{2}\\
    0 &  \zeta  &    z\\
    0 &    0    &    1
  \end{matrix}\right)
\mapsto
\left(
  \begin{matrix}
    \zeta  &    z\\
      0    &    1
  \end{matrix}\right).
\end{equation}
is a group homomorphism, which can be interpreted as a map
$$
  \Phi:Stab^{par}_G(x)\rightarrow \textrm{Isom}(\C)
$$
of the parabolic stabilizer onto the group of Euclidean
isometries of $\C$, with central kernel.

Since we assume that $\Gamma$ is a lattice, $Stab_\Gamma(x)$ must act
dicretely and cocompactly on the Heisenberg group $\Hc$, from which is
follows that $\Phi(Stab_\Gamma(x))$ must act dicretely and cocompactly
on $\C$ (and the kernel must be an infinite cyclic group).

In those terms, the analogue of the torsion-free condition for cusps
is quite natural; one requires that $\Phi(Stab_\Gamma(x))$ be a
torsion-free (cocompact) group of isometries of $\C$. It is a standard
fact that these are just free Abelian groups generated by two
translations in different directions, hence the quotient
$\C/\Phi(Stab_\Gamma(x))$ is an elliptic curve.

\begin{dfn}\label{dfn:neat}
  A lattice $\Gamma\subset PU(2,1)$ is {\bf neat} if every matrix
  representative $A$ of an element $\gamma\in\Gamma$ that has a root
  of unity as an eigenvalue is actually (a multiple of) a unipotent
  matrix.
\end{dfn}

One can push the above discussion to prove the following, see
~\cite{amrt} for the arithmetic case,~\cite{mok-projective} for the
general case.
\begin{thm}
  Let $\Gamma\subset PU(2,1)$ be a neat lattice. Then
  $X=\Gamma\backslash\CH^2$ is smooth, and it admits a smooth
  compactification $\bar X$, where $\bar X\setminus X$ is a disjoint
  union of elliptic curves with negative self-intersection.
\end{thm}

One can also give a simple interpretation of the self-intersection of
the elliptic curves in the compactification $\bar X$.

Indeed, let $\Gamma$ be a neat lattice and let $P$ be a cusp subgroup
of $\Gamma$. It follows from the previous discussion that $P$ is a
central extension (with infinite cyclic center) of sugroup of
$\textrm{Isom}(\C)$ generated by two translations. In terms of the
Heisenberg group, we get two Heisenberg translations $A,B$ (such that
the entries $A_{1,2}$ and $B_{1,2}$ linearly independent over $\R$)
such that $\Phi(A)$ and $\Phi(B)$ generate $\Phi(Stab_\Gamma(x))$,
and a non-trivial vertical translation $Z$ that generates the kernel.

Since the commutator of two Heisenberg translations must be a vertical
translation, we can write
$$
 ABA^{-1}B^{-1}=[A,B]=Z^k
$$
for some integer $k$. Proposition 4.2.12 and equation~(4.2.15)
of~\cite{holzapfel-book} say that the self-intersection of the
elliptic curve corresponding is $-|k|$.

This self-intersection can be computed very efficiently from a
presentation for $C$ in terms of generators and relations. Indeed,
given the above description, $C$ must be isomorphic to
$$
\FPG{a,b,z}{[a,b]z^{-k},\ z\textrm{ central }},
$$
whose abelianization is $\Z\oplus\Z\oplus\Z^{|k|}$.

\begin{rmk}
  When $\Gamma$ is a non-neat lattice, one can still compactify the
  quotient, but the ends are then filled with quotients of elliptic
  curves by a finite group, see recent work of
  Eyssidieux~\cite{eyssidieux}.
\end{rmk}

\subsection{Lattice triangle groups} \label{sec:triangle-groups}

Recall that complex reflections are elliptic elements whose matrix
representative has a repeated eigenvalue, see
equation~\eqref{eq:refl}. The lattices we consider in this paper are
all triangle groups in the following sense.
\begin{dfn}
  A complex hyperbolic triangle group is a subgroup of $PU(H)$
  generated by three complex reflections $R_1,R_2,R_3$.
\end{dfn}

It turns out that several triangle groups are lattices, as was first
observed by Mostow~\cite{mostow-pacific}. More examples were found
later~\cite{deraux-4445},~\cite{thompson-thesis},~\cite{dpp1}
and~\cite{dpp2}.

We refer to the lattice triangle groups listed in~\cite{dpp2} as
\emph{the known lattice triangle groups}. The list is reproduced here
in Tables~\ref{tab:neat-3s1} through~\ref{tab:neat-m3-0}.

We briefly review notation, a detailed parametrization for these
groups can be found in section~3 of ~\cite{dpp2}. We list groups of
three kinds, see Table~\ref{tab:triangles}.
\begin{table}[htbp]
  \begin{tabular}[htbp]{ccc}
    Name            &   Notation       & Parameters\\
    Sporadic groups & $\Sc(p,\tau)$    & $p\in\N^*,\tau\in\C$\\
    Thompson groups & $\Tc(p,{\bf T})$ & $p\in\N^*,{\bf T}\in\C^3$\\
    Mostow groups   & $\Gamma(p,t)$    & $p\in\N^*,t\in\Q$\\
  \end{tabular}
  \label{tab:triangles}
\end{table}
Sporadic groups and Mostow groups are generated by $R_1$ and $J$,
where $R_1$ is a complex reflection with multiplier
$e^{2\pi i/p}$, and $J$ is a regular elliptic element of order 3. One
gets two other complex reflections by setting $R_2=JR_1J^{-1}$ and
$R_3=J^{-1}R_1J$.

Sporadic groups are characterized by this data and the fact that
$\tr(R_1J)=\tau$. The corresponding sporadic group $\Sc(p,\tau)$ is a
lattice only for wisely chosen pairs $(p,\tau)$.

Mostow groups $\Gamma(p,t)$ are special cases of sporadic groups, where we take
$$
\tau=e^{\pi i(\frac{3}{2}+\frac{1}{3p}-\frac{t}{3})}.
$$

Thompson groups are generated by three complex reflections
$R_1,R_2,R_3$ with the same multiplier $e^{2\pi i/p}$, but where there
is no isometry $J$ as above such that $R_{k+1}=JR_kJ^{-1}$. These are
parametrized by a triplet ${\bf T}=(\rho,\sigma,\tau)\in\C^3$ that
generalizes the trace parameter of sporadic groups (see section~3
of~\cite{dpp2} for details).

As discussed in~\cite{dpp2}, the commensurability classes of the
lattices in these three classes contain all known non-arithmetic
commensurability classes of lattices in $PU(2,1)$ (there are currently
22 known commensurability classes).

\section{Conjugacy classes of isotropy groups} \label{sec:isotropy}

In this section we briefly recall how to use a fundamental domain $F$
for a lattice $\Gamma\subset PU(2,1)$ to obtain the list of conjugacy
classes of
\begin{itemize}
\item isotropy groups in $\Gamma$ of points in $\CH^2$;
\item cusps in $\Gamma$.
\end{itemize}
The isotropy groups of points in $\CH^2$ are precisely the maximal
finite subgroups in $\Gamma$, and the cusps are isotropy groups of
ideal boundary points that contain at least one parabolic element (the
latter are infinite, in fact they act cocompactly on horospheres).

Note that the finite isotropy groups were already listed in the tables
of~\cite{dpp2}; either they are generated by complex reflections (in
which case they occur in the tables giving vertex stabilizers), or
they are not (in which case they can be deduced from the information
about singular points of the quotient).

For cusps, the information given in~\cite{dpp2} only says whether the
group is generated by reflections, and when it is not, we did not give
actual generators. In order to get generators, we use the computer
output of Spocheck~\cite{spocheck}.

What we do in Spocheck uses the fundamental domains constructed
in~\cite{dpp2} and tracks cycles of (ideal or finite) vertices, as
well as facets of higher-dimension. The point is that every isotropy
group is conjugate to one with fixed point in the boundary of the
fundamental domain, hence it is enough to compute the stabilizers of
all facets of the fundamental domain. We use the same method as Mostow
in section 18.2 of~\cite{mostow-pacific}.

In the following paragraphs, we let $F$ be (finite-sided) a
fundamental domain for $\Gamma$, and let $v$ denote a facet of $F$. In
order to compute $Stab_\Gamma(v)$ we build a graph whose vertices are
given by facets of the fundamental domain $F$. For every such facet
$v$, we list all sides of $F$ containing $v$. For every such side $s$,
we consider the side-pairing map $\gamma_s$, and draw an edge from $v$
to $\gamma_s(v)$. Let us denote by $T$ the corresponding directed
graph; it is a finite graph since we assume $F$ has finitely many
sides.

By construction, the $\Gamma$-orbits of ideal vertices are in 1-1
correspondence with the connected components of $T$. We have a
well-defined representation $\rho_v:\pi_1(T,v)\rightarrow \Gamma$, and
the stabilizer of $v$ in $\Gamma$ is precisely the image
$Im(\rho_v)$. In particular, in order to get a generating set for
$Stab_\Gamma(v)$, it is enough to construct explicit generators of
$\pi_1(T,v)$ (which can be done by constructing a maximal subtree
containing $v$ in $T$).

The computations (and even the end results of these computations) are
too long to be included in a paper, we will list them in the form of a
Magma file in~\cite{computer-code}. In this paper, we give the details
only for a couple examples that illustrate the method, see
sections~\ref{sec:3s1} and~\ref{sec:3s5}.

The general result for ideal vertices is given in the Tables of
section~\ref{sec:neat} (Tables~\ref{tab:neat-3s1}
through~\ref{tab:neat-m3-0}), sixth (cusp generators) and seventh
column (cusp relations). 

\subsection{Isotropy groups for $\Gamma=\Sc(3,\sigma_1)$} \label{sec:3s1}

In this section, we use word notation used in~\cite{dpp2}, so that
$1$, $2$, $3$ and $4$ stand for $R_1$, $R_2$, $R_3$ and $J$
respectively, and $\bar 1$ stands for $R_1^{-1}$, etc.  Consider the
group $\Gamma=\Sc(3,\sigma_1)$. The Spocheck output~\cite{spocheck}
gives us the stabilizers listed in Table~\ref{tab:stabs-3s1}. It turns
out all non Abelian finite stabilizers are complex reflection groups,
so they can be described in terms of the Shephard-Todd
classification~\cite{shephard-todd}. In the tables, we write $G_k$ for
(a group isomorphic to) the $k$-th imprimitive
group in the Shephard-Todd list.\\

\begin{table}[htbp]
  \centering
  \begin{tabular}{c|c|c}
  Ridge stab.                    & Edge stab.                    & Vertex stab.\\
    \hline
  $\langle J\rangle$, order 3    & $\langle R_1\rangle$, order 3 & $\langle 1,2\rangle$, Cusp\\
  $\langle R_1\rangle$, order 3  &                               & $\langle 1,2\bar3\bar2\rangle$, order 72 ($G_5$)\\
  $\langle R_1J\rangle$, order 8 &                               & $\langle 1,232\bar3\bar2\rangle$, order 24 ($G_{4}$)\\
                                 &                               & $\langle 1,\bar3\bar2 323\rangle$, order 24 ($G_{4}$)
  \end{tabular}
  \caption{Facet stabilizers for $\Sc(3,\sigma_1)$}
  \label{tab:stabs-3s1}
\end{table}

Recall that $\br_n(a,b)$ means $(ab)^{n/2}=(ba)^{n/2}$, which when $n$
is odd means
\begin{equation}
  \label{eq:braiding}
  aba\cdots ba = bab\cdots ab 
\end{equation}
where both sides of equation~\eqref{eq:braiding} are products of $n$
factors.

For the cusp, we know from~\cite{dpp2} or~\cite{spocheck} that $R_1$
and $R_2$ braid with length 6. We then use the fact that $(ab)^3$ is
central in the braid group $\FPG{a,b}{\br_6(a,b)}$ and the basic
geometry of braiding complex reflections (see section~2.3
of~\cite{mostow-pacific} or section~2.3 of~\cite{dpp2}) to identify
the relevant subgroup of $\textrm{Isom}(\C)$ (see
section~\ref{sec:lattices}) as a $(3,3,3)$-triangle group.

From this and the fact that the braid relation $\br_6(a,b)$ implies
that $(ab)^3$ is central, one easily sees that the cusp has a
presentation of the form
$$
   \FPG{r_1,r_2}{r_1^6,r_2^6,\br_6(r_1,r_2)}.
$$

\subsection{Isotropy groups for $\Gamma=\Sc(3,\sigma_5)$} \label{sec:3s5}

A more complicated example is given by the group
$\Gamma=\Sc(3,\sigma_5)$. The Spocheck output gives us the information
in Table~\ref{tab:stabs-3s5}, where $\Xi$ is a complicated generating set for the corresponding
cusp. 
\begin{table}[htbp]
  \centering
  \begin{tabular}[htbp]{c|c|c}
     Ridge stab.                                             & Edge stab.                    & Vertex stab.\\
    \hline
    $\langle J\rangle$, order 3                              & $\langle R_1\rangle$, order 3 & $\langle 1,2\rangle$, order 72 ($G_5$)\\
    $\langle R_1\rangle$, order 3                            &                               & $\langle 1,2\bar3\bar2\rangle$, order 360 ($G_{20}$)\\
    $\langle R_1J\rangle$, order 30                          &                               & $\langle\Xi\rangle$, Cusp\\
    $\langle 2\bar3\bar2 123\bar2,(R_1J)^5\rangle$, order 18 &                               & \\
  \end{tabular}
  \caption{Facet stabilizers for $\Sc(3,\sigma_5)$}
  \label{tab:stabs-3s5}
\end{table}
Explicitly $\Xi$ consists of the following elements:
\begin{equation*}
\begin{array}{c}
  x_0=2\\
  x_1=123\bar2 12\bar3\,\bar2\,\bar1\\
  x_2=(123)^2 12\bar1 31\bar2\,\bar1(\bar3\,\bar2\,\bar1)^2\\
  x_3=(123)^4 1\bar323\bar 1(\bar3\,\bar2\,\bar1)^4\\
  x_4=(\bar3\,\bar2\,\bar1)^2\bar3\,\bar2\,31\bar3 23(123)^2\\
  x_5=(\bar3\,\bar2\,\bar1)23\bar2(123)\\
  x_6=\bar 4^2 1 2 1 \bar2 \, \bar1 \, \bar3 \, \bar2 \, \bar1 \, 2
\end{array}
\end{equation*}

By constructing all words of length $\leq 5$ in $x_0$ and $x_1$, one
checks that $x_2,x_3,x_4,x_5\in \langle x_0,x_1\rangle$, in fact we
have: $x_2=x_1x_0x_1^{-1}$, $x_3=x_1x_0x_1x_0^{-1}x_1^{-1}$,
$x_4=x_0^{-1}x_1^{-1}x_0x_1x_0$, $x_5=x_0^{-1}x_1x_0$.

Computing short words in $x_0,x_6$, we also have $x_1=x_6x_0x_6^{-1}$,
so we get the following.
\begin{prop}\label{prop:xi}
  The group generated by $\Xi$ is generated by $x_0$ and $x_6$.
\end{prop}

In order to work out a presentation for the group generated by $x_0$
and $x_6$, we need to work a little more. We use the discussion given
in section~\ref{sec:heisenberg} of the stabilizer of an ideal point in
$PU(2,1)$.

We start by giving matrices for the group, with integral entries in
the smallest possible number field, i.e. $\IK=\Q(\omega,\varphi)$,
where $\omega=\frac{-1+i\sqrt{3}}{2}$,
$\varphi=\frac{1+\sqrt{5}}{2}$. The lattice can be written as the
group generated by
\begin{equation*}
    \begin{array}[htbp]{c}
      R_1=\left(\begin{matrix}
          \omega & \omega+\varphi & \bar\omega\varphi+\omega\\
          0      &     1          &     0\\
          0      &     0          &     1
          \end{matrix}\right)
      J=\left(\begin{matrix}
          0      &     0          &   -\omega\\
          1      &     0          &     0\\
          0      &     1          &     0
          \end{matrix}\right),
    \end{array}
  \end{equation*}
which preserve the Hermitian form
$$
H=\left(\begin{matrix}
    \alpha            & \beta   & -\omega\bar\beta\\
    \bar\beta         & \alpha  &  \beta\\
    -\bar\omega\beta  & \bar\beta  & \alpha
  \end{matrix}\right),
$$
where $\alpha=3$ and $\beta=-(2+\omega)\varphi+(1-\omega)$.

It follows from Proposition~\ref{prop:xi} that that the cusp is
generated by
$$
A = \bar 2\,12312\bar 1\,\bar 2\,\bar 1\,P^2, B = R_2
$$
By using a suitable matrix $Q\in GL(3,\Oc_\IK)$, for example
$$
Q = \left(\begin{matrix}
    \omega(\varphi+1)       & \bar\omega(\varphi+1)      & \omega+1\\
    -\bar\omega(2\varphi+1) & -(\varphi+1)               & \varphi+\bar\omega\\
    \varphi+1               &  \omega(\varphi+1)         & 1+\bar\omega\varphi
  \end{matrix}\right),
$$
we can write $A$ and $B$ in upper triangular form, namely we have
(in the projective group)
$$
Q^{-1}AQ=\left(
  \begin{matrix}
    1 &        0    &  -\bar\omega\\
    0 & -\bar\omega &  0\\
    0 &    0        &  1
  \end{matrix}
\right),\quad 
Q^{-1}BQ=\left(
  \begin{matrix}
    1 & 1+2\omega & -\omega\\
    0 & \omega   &  1\\
    0 & 0        &  1
      \end{matrix}
\right). 
$$
Note that $Q^*HQ$ is (a real multiple) of
$$
\left(
  \begin{matrix}
    0            &  0  & \omega-1\\
    0            &  3  &  0\\
    \bar\omega-1 & 0 & 0
  \end{matrix}
\right).
$$

As discussed in section~\ref{sec:heisenberg}, the action on the complex
line tangent to the ideal boundary at the ideal fixed point is
described projectively by the lower-right $2\times 2$ submatrix,
i.e. if we denote by $\tilde{A},\tilde{B}$ the corresponding Euclidean
isometries of $\C$, we have
$$
\tilde{A}(z)=-\bar\omega z,\quad \tilde{B}(z)=\omega z+1
$$
which are rotations about $0$ (resp. $\frac{1-\bar\omega}{3}$) and
order $6$ (resp. $3$).

Note that the product $\tilde{B}\tilde{A}$ is a rotation of order
2, and the relevant Euclidean group of isometries is a $(2,3,6)$
triangle group.

Going back to the original matrices $A,B$ and $BA$, note that
(still in the projective group)
\begin{equation}\label{eq:powers}
A^6=\left(
  \begin{matrix}
    1 & 0 & -6\bar\omega\\
    0 & 1 & 0\\
    0 & 0 & 1
  \end{matrix}\right),\quad
  B^3=Id,\quad
  (BA)^2 = 
  \left(\begin{matrix}
    1 & 0 & -\bar\omega\\
    0 & 1 & 0\\
    0 & 0 & 1
  \end{matrix}\right).
\end{equation}
The last equality shows that $(BA)^2$ is central.

In fact we have the following
\begin{prop}
  The cusp of $\Sc(3,\sigma_5)$ can be represented as the group
  generated by $A$ and $B$. Its center is generated by $(AB)^2$, and
  it is isomorphic to 
  $$
  \langle a,b | a^6(ba)^{-12}, b^3, [a,(ab)^2], [b,(ab)^2]\rangle.
  $$
\end{prop}
\begin{pf}  
  The above discussion shows that the group of isometries of $\C$
  obtained as the action on the complex line tangent to the ball at
  the fixed point of the cusp is a $(2,3,6)$-triangle group.

  It is well known that this triangle group has a presentation of the
  form
  $\FPG{\alpha,\beta}{\alpha^6,\beta^3,(\alpha\beta)^2}$,
  in particular, any word in $\alpha,\beta,\gamma$ that is trivial in
  the triangle group can be written as a product of conjugates of
  $\alpha^6$, $\beta^3$ or
  $(\alpha\beta)^2$. This implies that every central element in the cusp must
  be a product of powers of $A^6$, $B^3$ and
  $(AB)^2$. The computation of equation~\eqref{eq:powers} shows that
  $(AB)^2$ generates the center.
\end{pf}

Note that explicit computation shows that
$BA^{-1}$ has order has order 6, and that the braid relation
$\br(B,BA^{-1})$ holds in the group. One checks (for example using
Magma) that the group
$$
 \FPG{ c,d }{ c^3,d^6,\br_4(c,d) }
$$
is isomorphic to the above presentation (where the isomorphism
maps $c\leftrightarrow b$, $d\leftrightarrow ba^{-1}$).

\section{Methods for finding subgroups} \label{sec:methods}

In order to get subgroups, we use four basic methods, listed
\begin{itemize}
\item The \verb|LowIndexSubgroups| gives us list of subgroups of small
  index (usually reasonable for index bounds of about 30);
\item The \verb|LowIndexNormalSubgroups| gives us list of \emph{normal}
  subgroups of small index, and it works for much larger index than
  the previous one (usually reasonable for index bounds of 30,000 to
  100,000);
\item The \verb|SimpleQuotients| allows us to study only simple
  quotients, and can give homomorphisms to very large finite (simple)
  groups;
\item Congruence subgroups (either we choose an explicit prime ideal
  and reduce the matrices modulo that ideal, or we use the
  \verb|CongruenceImage| command in Magma).
\end{itemize}
The cost of all these methods grows exponentially as the number of
generators increase. What makes them usable in the context of
subgroups of small index in this context is that lattice triangle
groups have few generators (most are generated by 2 elements, some by
3 elements).

There is also a slightly more subtle method to build non-normal (neat)
subgroups containing a given normal (neat) subgroup. This is well
known to experts, but we give some details in
section~\ref{sec:promotion}.

\subsection{Low Index Subgroups (LIS)}

This method works only for subgroups of fairly small index; the
indices accessible by this method depend to a great extend on the
lattice triangle group we are considering (generally index 20 is
already a lot to ask).

\subsection{Low Index Normal Subgroups (LINS)}

For normal subgroups, the allowed indices are much larger, but once
again, the reasonable values depend on the group (for some, Magma
gives reasonable access to the index $\leq 100\,000$, while for some
others index $1\,000$ is already a lot to ask).

\subsection{Simple Quotients (SQ)}

The search for simple quotients depends on choices of parameters
(lower and upper bounds for the order of the group, and an upper bound
on the degree of the permutation group, i.e. $N$ such that the group
embeds in the symmetric group $S_N$). For more details, see the README
file in the computer code~\cite{computer-code}.

\subsection{Congruence Subgroups (CS)}

We say a word about how we compute congruence subgroups. Each triangle
group comes with an explicit generating set, given by a finite set
$A$ of matrices in $U(H)$, with algebraic integer entries. Denote by
$\widetilde{\Gamma}$ the group generated by $A$, and by $\Gamma$ the
corresponding subgroup of $PU(H)$.

Here $H$ is a Hermitian form over a number field $\IK$, with ring of
integers $\Oc_\IK$. Given an ideal $I\subset \Oc_\IK$, we consider the
finite ring $R=\Oc_\IK/I$, and consider the group homomorphism
$$
  \varphi:\widetilde{\Gamma}\rightarrow GL(3,R)
$$
obtained by reducing all entries modulo $I$. 
\begin{dfn}\label{dfn:congruence-subgroup}
  The principal congruence subgroup of $\Gamma$ mod $I$ is the
  projectivization of the kernel of $\varphi$.
\end{dfn}

In order for these principal congruence subgroups to be accessible to
computation, we need to define the above data in Magma, which is a
little subtle.

Recall that we are given a presentation for $\Gamma$ in terms of
generators and relations, and we would like to compute a presentation
for a principal congruence subgroup. In order to do this, we follow
the next steps.
\begin{itemize}
\item Define $\widetilde{F}=\textrm{Im}(\varphi)$ as the
  \verb|MatrixGroup| over the finite ring $R=\Oc_\IK/I$ generated by
  $\varphi(a)$, $a\in A$.
\item Convert $\widetilde{F}$ to a permutation group via
  \verb|PermutationGroup(FPGroup(| \dots \verb|));|, and compute the
  permutation group $F=\widetilde{F}/Z_{\widetilde{F}}$, where
  $Z_{\widetilde{F}}$ is the group of scalar matrices in
  $\widetilde{F}$.
\item Define a homomorphism from $\Gamma$ to $F$, and compute its
  kernel.
\end{itemize}
The above steps are quite heavy computationally, and they will only
succeed when the order of $\widetilde{F}$ is not too large.

Note that in order to use the above method, we need to choose an ideal
$I$ in $\Oc_\IK$. Recall that
$\Oc_\IK$ is not a unique factorization domain, but there is a unique
factorization of ideals as a product of prime ideals, see standard
texts on algebraic number theory, e.g.~\cite{neukirch}. This
factorization is implemented in Magma.

In order to select ideals, we factor a the ideal
$n\Oc_\IK=I_1\cdot I_2\cdots I_r$ for some rational integer $n\in\Z$, and pick
one of its prime ideal factors $I_1$ (or sometimes the product of
several such prime ideal factors).

\subsection{Promoting normal subgroups to non-normal subgroups of smaller index}\label{sec:promotion}

The last three methods (LINS, QS, SC) give normal subgroups of fairly
large index, and we now explain how to improve this to get
(non-normal) subgroups of smaller index.

Suppose $\varphi:\Gamma\rightarrow F$ is a surjective morphism onto a
finite group, let $K=\textrm{Ker}(h)$. The basic observation is that
for every subgroup $S\subset F$, $H=\varphi^{-1}(S)$ is a subgroup of
$\Gamma$ that contains $K$. Converserly, any subgroup $H$ with
$K\subset H\subset \Gamma$ is obtained in this way (as
$\varphi^{-1}(\varphi(H))$). Note also that with the above notation,
the indices $[\Gamma:H]$ and $[F:S]$ are equal.

Moreover, the preimage $\varphi^{-1}(S)$ is a finitely presented group
and, at least when the index is not too large, a presentation can be
obtained via Magma.

Now suppose that we have a list $I_1,\dots,I_r$ of finite subgroups of
$\Gamma$ such that every torsion element in $\Gamma$ is conjugate to
an element of some $I_k$ (in other words, the subgroups $I_k$ give a
list of the non-trivial isotropy groups for the action of $\Gamma$ on
$\CH^2$).
\begin{prop}\label{prop:promote-tf}
  \begin{enumerate}
  \item $K$ is torsion-free if and only if $|\varphi(I_j)|=|I_j|$ for
    all $j=1,\dots,r$.
  \item If $K$ is torsion-free, then $H=\varphi^{-1}(S)$ is
    torsion-free if and only if $x^{-1}Sx\cap I_j=\{Id\}$ for every
    $x\in F$ and every $j$.
  \end{enumerate}
\end{prop}
In part~(2), instead of taking all the elements $x\in F$, it is of
course enough to check all the elements in a right-transversal for $S$
in $F$. Note that all the verifications of
Proposition~\ref{prop:promote-tf} take place in a finite group, so
in a sense they can be considered easy\dots but they can take a long time if the index
$[F:S]$ is big.

There is also a variant of this method that allows us to check
whether $K$ and $H$ are neat (see Definition~\ref{dfn:neat}), but this
is quite a bit more complicated. We now summarize its main steps.
\begin{itemize}
\item Let $C$ denote a cusp of $G$ (we discussed how to get generators
  and relations for $C$ in section~\ref{sec:3s5}).
\item Let $C_F$ denote $h(C)$ (generators for this group are given by
  the $h$-images of generators of $C$), and
  $K_C=\textrm{Ker}(h|C)$. Note that this has finite index in $C$, so
  Magma can find a presentation for $K_C$. Let $k_1,\dots,k_n$ denote
  a generating set for $K_C$.
\item The cusps of $K\backslash\CH^2$ are in 1-1 correspondence with right
  cosets of $C_F$ in $F$, and the corresponding cusp groups are all
  conjugate to $K_C$. Denote by $f_1,\dots,f_m$ a right transversal
  for $C_F$.
\item Study the action of $S=h(H)$ on the set of the right coset of
  $C_F$, find its orbits and stabilizers. The cusps of $H$ are in 1-1
  correspondence with these orbits.
\item In order to get cusp generators, take a right coset $x=C_Ff_j$
  and denote by $I_x$ its stabilizer under the $S$-action. Find a
  generating set for $I_x$, and choose lifts $i_1,\dots,i_r$ to $H$
  for this generating set. Cusp generators are given by the elements
  $$
  i_1,\dots,i_r,k_1^{f_j},\dots,k_n^{f_j}.
  $$
\item Use a suitable linear change of coordinates to make all cusp
  generators upper triangular, and multiply each of them by a scalar
  to get the upper left entry to be 1 (recall we are working in
  $PU(2,1)$ rather than $U(2,1)$). The corresponding cusp is neat if
  and only if every generator is unipotent.
\item Since every cusp is isomorphic to a finite index subgroup of
  $K_C$, we can use Magma to get a presentation for every cusp group.
\item If the cusp is neat, we compute the self-intersection of
  corresponding the elliptic curve in the toroidal compactification by
  computing the abelianization of the corresponding cusp group, which
  is isomorphic to $\Z^2\oplus\Z_q$ for some $q\in \N^*$. In that
  case the self-intersection is given by $-q$ (see Proposition~4.2.12
  in~\cite{holzapfel-book}).
\end{itemize}

A Magma implementation of these methods can be found in our computer
code~\cite{computer-code}.

\section{A tower of one-cusped smooth ball
  quotients}\label{sec:one-cusp}

In this section we describe prove Theorem~\ref{thm:one-cusp}. As
mentioned in the introduction, this gives a positive answer to a
question raised by Stover in~\cite{stover-cusps}.

The group is a subgroup of index 72 in the Mostow group
$G=\Gamma(6,0)$, which has the presentation
$$
\FPG{r_1,r_2,r_3,j}{r_2^{-1}jr_1j^{-1},r_3^{-1}j^{-1}r_1j, r_1^6, j^3, (r_1j)^{12}, (r_2r_1j)^6, \br_3(r_1,r_2)}.
$$
This group has a single cusp, represented by $C=\langle r_1,r_2\rangle$.

Every non-trivial isotropy group for the action of $G$ on $\CH^2$ is
in the list of table~\ref{tab:isotropy}.
\begin{table}[htbp]
  \centering
  \begin{tabular}[htbp]{c||c|c|c|c|c|c}
    Isotropy group & $\langle j\rangle$ & $\langle r_1j\rangle$ & $\langle r_1,r_2r_1j\rangle$ & $\langle r_2j^{-1}\rangle$ & $\langle r_2,(r_1j)^2\rangle$ & $\langle r_2r_1j, (r_1j)^2\rangle$\\
    \hline
    Order &       $3$         &         $12$          &              $36$             &           $12$           &         $36$                  &       $36$
  \end{tabular}
  \caption{Isotropy groups of $G$}
  \label{tab:isotropy}
\end{table}

The subgroup alluded to in the statement of Theorem~\ref{thm:one-cusp}
is
$$
H=\langle r_1r_2r_1jr_2^{-1}, r_1r_2r_3jr_1^{-1}, jr_1^{-1}r_3r_2r_3, (r_1^2r_2)^2\rangle.
$$
Alternatively, the reader can reconstruct the group by loading the
code given in~\cite{computer-code} into Magma, then running the
commands
\begin{center}
\begin{tabular}{l}
  \verb|A:=MostowGroup(6,0);|\\
  \tt{FindTorsionFreeSubgroups($\sim$A,864,[2]);}
\end{tabular}
\end{center}
The field \verb|A`Subgroups| should then have a seven elements, and
\verb|A`Subgroups[7]| contains a Magma description of the subgroup in
question.
We write $X=H\backslash \CH^2$.

From the above pieces of information, Magma allows us to check that
$[\Gamma:H]=72$, hence 
$$
\chi(X)=72\cdot\chi^{orb}(\Gamma\backslash\CH^2)=72\cdot\frac{1}{12}=6.
$$

Let $K=\textrm{Core}_G(H)$, and let $F=\Gamma/K$, which has order
$864$. Using Magma, it is straighforward to check that the isotropy
groups of Table~\ref{tab:isotropy} inject in $F$, so $K$ is
torsion-free.

Let us denote by $S\subset F$ the image of $H$. A simple computer
check shows that $S$ has trivial intersection with every conjugate of
the image in $F$ of every isotropy group in Table~\ref{tab:isotropy},
so $H$ is torsion-free.

We write $X=H\backslash \CH^2$ and $Y=K\backslash\CH^2$, and $\bar X$,
$\bar Y$ for their respective toroidal compactifications.

In order to study the cusps of $X$ and $Y$, we let $C_F$ denote the
image in $F$ of the cusp $\langle r_1,r_2\rangle$ in $G$; we have
$[F:C_F]=12$, so $Y$ has 12 cusps.  Computing the kernel $CK$ of the
homomorphism $h_C:C\rightarrow C_F$ using Magma, one computes
generators, and checks that each generator is unipotent element, so
$K$ is neat. The Abelianization of $CK$ is $\Z_{12}\oplus\Z^2$, so the
$12$ elliptic curves compactifying $Y$ to $\bar Y$ have
self-intersection $-12$.

In order to prove that $X$ has a single cusp, simply check that $S$
acts transitively on the set of right cosets of $C_F$ in $F$. It also
follows that the 12 cusp groups of $X$ are all isomorphic to the cusp
group $CK$ of $X$. Using Magma, one can check that the cusp $C$ of $H$
is represented by the group generated by
\begin{equation}\label{eq:cusp-gens}
  (r_2r_1^2)^{2},\quad r_2^2r_1^{-1}r_2r_1^{-2}r_2r_1^{-1},\quad r_2r_1^{-2}r_2r_1^{-1}r_2^2r_1^{-1}.
\end{equation}

We finish this section by proving Theorem~\ref{thm:tower}. Let
$\varphi:H\rightarrow \Z^2$ be obtained from the abelianization map by
projecting onto the $\Z^2$-factor of $\Z_3\oplus\Z^2$. By explicit
calculation of the image of the generating set for $C$ of
equation~\eqref{eq:cusp-gens}, we check that $\varphi(C)$ is a lattice
in of index $12$ in $\Z^2$.

By projecting onto one factor of $\Z^2$, it is easy to get a
homomorphism $\psi:H\rightarrow\Z$ such that $\psi(C)=2\Z$. Now for
$n\in\N^*$, consider $H_n=\psi^{-1}(n\Z)$. For any odd $n$, the action
of $n\Z$ is transitive on the cosets of $2\Z$, so $H_n$ has exactly
one cusp. The tower of Theorem~\ref{thm:tower} is then obtained by
choosing an increasing sequence of odd integers $n_1<n_2<\dots$ such
that for all $j$, $n_j$ divides $n_{j+1}$.

\section{Summary of the results for all groups} \label{sec:results}

In this section, we list the subgroups we were able to find inside the
known lattice triangle groups. We focus on three different aspects,
namely
\begin{itemize}
\item neat subgroups (section~\ref{sec:neat})
\item subgroups with positive $b_1$ (section~\ref{sec:betti})
\item large subgroups (section~\ref{sec:large})
\end{itemize}

\subsection{Positive $b_1$} \label{sec:betti}

In order to find subgroups $\Gamma_0\subset G$ with $b_1(\Gamma_0)>0$,
we study whichever neat subgroups we can find in $G$, ask Magma to
compute the Abelizanition using the command \verb|AbelianQuotientInvariants|.

When getting a homomorphism by using \verb|SimpleQuotients| or
\verb|CongruenceImage|, the large order of the group often seems to
make it too difficult for Magma to compute a presentation for
$K=Ker(h)$, it is not clear how to estimate $b_1(K)$.

In Tables~\ref{tab:betti-sporadic} through~\ref{tab:betti-mostow}, we
list the groups where one of these methods yields an explicit
(torsion-free or not) subgroup with positive $b_1$. In case we have
such an example we mention how the reader can find it, writing
\begin{itemize}
\item "LIS(N)" where $N$ is the index of the subgroup;
\item "LINS(N)" where $N$ is the (bound on the) index of the normal
  subgroup or "LINS(N,G)" with an concise description of the quotient
  if we have one;
\item "SQ(G)" where $G$ is a description of the Simple Quotient;
\item "mod $p$" if reduction mod $p$ gives a kernel with positive $b_1$.
\item ``Map(N)'' if we know a homomorphism to a triangle group having
  a torsion-free subgroup of index $N$ (see section~\ref{sec:large}).
\item nothing if we were unable to find any subgroup with positive $b_1$.
\end{itemize}

For non-cocompact groups, in the last column of the tables, we give
information as to whether or not the positive first Betti number comes
from the cusps; more precisely, we mention whether (for a suitable
subgroup in the list) all cusp groups map have infinite image in the
abelianization.

\begin{table}[htbp]
  \centering
  \begin{tabular}[htbp]{|c|c|c|c|c|c|}
    \hline
    NA & NC & $\Sc(3,\sigma_1)$     & OK & LINS(18144), SQ                                     & ?\\
    NA & NC & $\Sc(4,\sigma_1)$     & OK & LINS(96)                                            & $\infty$ cusp images\\
    NA & NC & $\Sc(6,\sigma_1)$     & OK & LINS(6)                                             & $\infty$ cusp images\\
    \hline
       &    & $\Sc(3,\bar\sigma_4)$  & OK & LINS(18144)                                        & \\
    NA & NC & $\Sc(4,\bar\sigma_4)$  & ?  &                                                    & ? \\
    NA &    & $\Sc(5,\bar\sigma_4)$  & ?  &                                                    & \\
    NA & NC & $\Sc(6,\bar\sigma_4)$  & ?  &                                                    & ? \\
    NA &    & $\Sc(8,\bar\sigma_4)$  & ?  &                                                    & \\
    NA &    & $\Sc(12,\bar\sigma_4)$ & ?  &                                                    & \\
    \hline
       &    & $\Sc(2,\sigma_5)$  & ?  &                                                    & \\
    NA & NC & $\Sc(3,\sigma_5)$  & OK & LINS(360)                                          & $\infty$ cusp images\\
    NA & NC & $\Sc(4,\sigma_5)$  & ?  &                                                    & ? \\
    \hline
       &    & $\Sc(3,\sigma_{10})$   & OK  & LINS(2160)                                    & \\
       &    & $\Sc(4,\sigma_{10})$   & ?   &                                               & \\
       &    & $\Sc(5,\sigma_{10})$   & OK  & LINS(600)                                     & \\
       &    & $\Sc(10,\sigma_{10})$  & OK  & LINS(18000)                                   & \\
    \hline
  \end{tabular}
  \caption{Subgroups of Sporadic groups with $b_1>0$}
  \label{tab:betti-sporadic}
\end{table}

\begin{table}[htbp]
  \centering
  \begin{tabular}[htbp]{|c|c|c|c|c|c|}
    \hline
       &    & $\Tc(3,{\bf S}_2)$    & OK & LINS(360)       & \\
    NA & NC & $\Tc(4,{\bf S}_2)$    & OK & LINS(23040)     & not all cusps have $\infty$ image\\
       &    & $\Tc(5,{\bf S}_2)$    & ?  &                 & \\
    \hline
       & NC & $\Tc(3,{\bf E}_2)$    & OK & LINS(24)        & $\infty$ cusp images \\
    NA & NC & $\Tc(4,{\bf E}_2)$    & OK & LINS(24)        & $\infty$ cusp images\\
       & NC & $\Tc(6,{\bf E}_2)$    & OK & LINS(6)         & $\infty$ cusp images \\
    \hline
       &    & $\Tc(2,{\bf H}_1)$    & OK & LIS(56)         & \\
    \hline
       &    & $\Tc(2,{\bf H}_2)$    & OK & LIS(30)         & \\
       &    & $\Tc(3,{\bf H}_2)$    & ?  &                 & \\
       &    & $\Tc(5,{\bf H}_2)$    & OK & LINS(600)       & \\
    \hline
  \end{tabular}
  \caption{Subgroups of Thompson groups with $b_1>0$}
  \label{tab:betti-thompson}
\end{table}

\begin{table}[htbp]
  \centering
  \begin{tabular}[htbp]{|c|c|c|c|c|c|}
    \hline
       &    & $\Gamma(5,7/10)$   & OK & LINS(15000) & \\
       & NC & $\Gamma(6,2/3)$    & OK & LINS(6)     & $\infty$ cusp images\\
       &    & $\Gamma(7,9/14)$   & OK & Map(84)     & \\
       &    & $\Gamma(8,5/8)$    & OK & LIS(12), Map(48) & \\
       &    & $\Gamma(9,11/18)$  & OK & LIS(9)      & \\
       &    & $\Gamma(10,3/5)$   & OK & LIS(10), Map(30) & \\
       &    & $\Gamma(12,7/12)$  & OK & LIS(6), Map(72) & \\
       &    & $\Gamma(18,5/9)$   & OK & LIS(6), Map(54) & \\

       &    & $\Gamma(4,5/12)$   & OK & LINS(18144) & \\
       &    & $\Gamma(5,11/30)$  & ?  &             & \\
       & NC & $\Gamma(6,1/3)$    & OK & LINS(6)     & $\infty$ cusp images\\
       &    & $\Gamma(7,13/42)$  & OK & Map(84)     & \\
       &    & $\Gamma(8,7/24)$   & OK & LIS(12), Map(48) & \\
       &    & $\Gamma(9,5/18)$   & OK & LIS(9)      & \\
       &    & $\Gamma(10,4/15)$  & OK & LIS(10), Map(30) & \\
       &    & $\Gamma(12,1/4)$   & OK & LIS(6), Map(72) & \\
       &    & $\Gamma(18,2/9)$   & OK & LIS(6), Map(18) & \\
    
       &    & $\Gamma(3,1/3)$    & OK & LINS(18144) & \\
       & NC & $\Gamma(4,1/4)$    & OK & LINS(96)    & $\infty$ cusp images\\
       &    & $\Gamma(5,1/5)$    & ?  &             & \\
    NA & NC & $\Gamma(6,1/6)$    & OK & LIS(6), Map(18) & $\infty$ cusp images\\
       &    & $\Gamma(8,1/8)$    & OK & LINS(1536), Map(72) & \\
       &    & $\Gamma(12,1/12)$  & OK & LIS(6), Map(72) & \\

       &    & $\Gamma(3,7/30)$   & ?  &             & \\
       &    & $\Gamma(4,3/20)$   & OK & Map(120)    & \\
       &    & $\Gamma(5,1/10)$   & OK & LINS(600)   & \\
       &    & $\Gamma(10,0)$     & OK & LINS(600), Map(30) & \\

       & NC & $\Gamma(3,1/6)$    & OK & LINS(72)    & $\infty$ cusp images\\
       &    & $\Gamma(4,1/12)$   & OK & LINS(864), Map(18) & \\
       & NC & $\Gamma(6,0)$      & OK & LINS(6)     & $\infty$ cusp images\\

       &    & $\Gamma(3,5/42)$   & ?  &             & \\

       &    & $\Gamma(3,1/12)$   & OK & LIS(48)     & \\
       &    & $\Gamma(4,0)$      & OK & LIS(24), Map(48) & \\

       &    & $\Gamma(3,1/18)$   & OK & LIS(24)     &\\

       &    & $\Gamma(3,1/30)$   & OK & LIS(40)     & \\

       &    & $\Gamma(3,0)$      & OK & LIS(24)     & \\
    \hline
  \end{tabular}
  \caption{Subgroups of Mostow groups with $b_1>0$}
  \label{tab:betti-mostow}
\end{table}

\subsection{Large subgroups}\label{sec:large}

We call a group $G$ large if it has a finite index subgroup $H$ with a
surjective homomorphism $\varphi:H\rightarrow F_n$, where $n>1$ and
$F_n$ is the non-Abelian free group on $n$ generators. Of course, we
may assume $n=2$, since every $F_n,n>2$ maps onto $F_2$. If this is
the case, then by abelianizing the free group, we get a morphism onto
$\Z^2$, hence the first Betti number of $H$ is positive.

Even though there are nice sufficient conditions of largeness due to
Button~\cite{button}, there is no general algorithm for determining
largeness. The sufficient condition for largeness used in this paper
is based on the following simple observation.
\begin{lem}\label{lem:sl2r}
  \begin{enumerate}
  \item Every torsion-free lattice in $PSL(2,\R)$ maps onto a
    non-Abelian free group;
  \item Every lattice in $PSL(2,\R)$ is large.
  \end{enumerate}
\end{lem}
\begin{pf}
  Part~(2) follows from part~(1) because of Selberg's lemma. Part~(1)
  is obvious for non cocompact lattices, since the fundamental group
  of a non-compact hyperbolic Riemann surface is itself a non-Abelian
  free group. For cocompact lattices, we can write the corresponding
  closed surface as the boundary of a handlebody, which retracts onto
  a bouquet of circles.
\end{pf}
As an immediate consequence, we see that if $G$ maps onto a lattice in
$PSL(2,\R)$, then $G$ is large.

Note that $PSL(2,\R)\cong PU(1,1)$, but when we change this group to
$PU(2,1)$, we have no reasonable analogue of Lemma~\ref{lem:sl2r} at
hand. It is a folklore conjecture that every lattice in $PU(2,1)$ of
Kazhdan type (i.e. such that its ambient algebraic $\Q$-group is the
group of a Hermitian form over a number field, rather than a more
complicated division algebra) should have a large subgroup of finite
index. Every complex reflection group is of Kazhdan type, so all the
groups we consider in this paper fall in the scope of this conjecture,
and we expect them to be large; but for most known lattice triangle
groups, largeness is not known.

We now discuss how to use Lemma~\ref{lem:sl2r} in order to find
explicit subgroups of finite index in some Mostow groups that map to
non-Abelian free groups. We do not know how to generalize this to
other lattice triangle groups.

It is well known to experts that Livne's thesis~\cite{livne} gives
maps from \emph{some} Mostow lattices to surface groups, which shows
that some specific Mostow lattices are large. For a description of
Livne's lattices in relation to Mostow's groups, see~\S 16 in the book
by Deligne and Mostow~\cite{deligne-mostow-book}. Recall that Livne
constructs smooth ball quotients $X_n$ indexed by
$n=5,6,7,8,9,10,12,18$, whose automorphism group $A$ fits is an
extension
$$
1\rightarrow \Z/d\Z \rightarrow A \rightarrow SL(2,\Z/n\Z)\times (\Z/n\Z)^2 \rightarrow 1,
$$
and the quotient $A\backslash X_n$ is
$\Gamma_{\mu,\Sigma}\backslash\CH^2$, where $d=gcd(6,n)$. The order of
these automorphism groups gives the index of the corresponding
subgroup of the relevant Mostow group; we list these in the last
column of Table~\ref{tab:livne}.
\begin{table}[htbp]
  \centering
  \begin{tabular}[htbp]{|c|c|c|c|c|}
    \hline
    Mostow group              & $n$  & $d$ & $|SL(2,n)|$ &  Index\\
    \hline
    $\Gamma(5,\frac{7}{10})$  & $5$  & $5$ &    120      &  15\,000\\
    $\Gamma(6,\frac{2}{3})$   & $6$  & $1$ &    144      &  5\,184\\
    $\Gamma(7,\frac{9}{14})$  & $7$  & $7$ &    336      &  115\,248\\
    $\Gamma(8,\frac{5}{8})$   & $8$  & $4$ &    384      &  98\,304\\
    $\Gamma(9,\frac{11}{18})$ & $9$  & $3$ &    648      &  157\,464\\
    $\Gamma(10,\frac{3}{5})$  & $10$ & $5$ &    720      &  360\,000\\
    $\Gamma(12,\frac{7}{12})$ & $12$ & $2$ &    1\,152   &  331\,776\\
    $\Gamma(18,\frac{5}{9})$  & $18$ & $3$ &    3\,888   &  3\,779\,136\\
    \hline
  \end{tabular}
  \caption{Index of torsion-free subgroups coming from Livne's construction}
  \label{tab:livne}
\end{table}
For some Mostow groups, we will give subgroups of much smaller index
that map onto a non-Abelian free group, see
Proposition~\ref{prop:maps-to-curves} and
Table~\ref{tab:triangle-groups}.

Note that the maps from Livne groups to Fuchsian groups are actually
induced by holomorphic fibrations of the corresponding ball quotients
over suitable Riemann surfaces (actually $X_n$ fibers over the
quotient of the upper half plane under the principal congruence
subgroup modulo $n$, see~\cite{deligne-mostow-book}). The maps we
construct are also induced by holomorphic maps to curves, as can be
seen from their interpretation as forgetful maps of moduli spaces of
points on $\CP^1$, see~\cite{deraux-maps3d}.

Rather than going into the details of Deligne-Mostow theory and
forgetful maps, we give a more down to earth approach that requires
only group theory.

\begin{prop}\label{prop:maps-to-curves}
  The Mostow groups $\Gamma(7,9/14)$, $\Gamma(8,5/8)$,
  $\Gamma(9,11/18)$, $\Gamma(10,3/5)$, $\Gamma(12,7/12)$,
  $\Gamma(18,5/9)$, $\Gamma(7,13/42)$, $\Gamma(8,7/24)$,
  $\Gamma(9,5/18)$, $\Gamma(10,4/15)$, $\Gamma(12,1/4)$,
  $\Gamma(18,2/9)$, $\Gamma(6,1/6)$, $\Gamma(8,1/8)$,
  $\Gamma(12,1/12)$, $\Gamma(4,3/20)$, $\Gamma(10,0)$,
  $\Gamma(4,1/12)$, $\Gamma(6,0)$, $\Gamma(4,0)$ have a surjective
  homomorphism onto a lattice in $PSL(2,\R)$.
\end{prop}

The basis of the construction of these homorphisms is simply to
combine Mostow's braid group description of the lattices $\Gamma(p,t)$
with the results in~\cite{deraux-maps3d} about forgetful maps of
Deligne-Mostow moduli spaces.

We briefly sketch how this is done; the reader without any knowledge
of Deligne-Mostow theory can skip this part. Each group
$\Gamma(p,t)=\langle R_1,J\rangle$ is conjugate in $PU(2,1)$ to a
specific Deligne-Mostow group
from~\cite{deligne-mostow},~\cite{mostow-ihes}, namely
$\Gamma_{\mu,\Sigma}$ where
$$
\mu=\left(\frac{1}{2}-\frac{1}{p},\frac{1}{2}-\frac{1}{p},\frac{1}{2}-\frac{1}{p},
  \frac{1}{4}+\frac{3}{2p}-\frac{t}{2},\frac{1}{4}+\frac{3}{2p}+\frac{t}{2}\right)
$$
and $\Sigma=S_3$ permutes the first three weights.

If $\mu$ satisfies the Deligne-Mostow INT condition, then the
corresponding Deligne-Mostow group $\Gamma_{\mu}$ is contained in
$\Gamma_{\mu,\Sigma}$ as a subgroups of index $|\Sigma|=3!=6$. Mostow
explains in~\cite{mostow-survey} how to write explicit generators for
$\Gamma_{\mu}$, namely $A_j,B_j$, $j=1,2,3$ with 
\begin{equation}
  \label{eq:generators-Gamma_mu}
  B_j=R_j^2,\quad A_j^{-1}=J^{-1}R_jR_{j+1}.
\end{equation}
Recall that $R_{j+1}=JR_jJ^{-1}$ index $j$ taken mod 3
($J^3=Id$). From this information, it is easy to ask Magma for a
presentation for each $\Gamma_\mu$.

Note also that when $\mu$ satisfies Mostow's $\Sigma$-INT condition
for $\Sigma=S_3$ but not the INT condition, we have
$\Gamma_\mu=\Gamma_{\mu,\Sigma}$, see~\cite{sauter} for instance.

Recall the presentation for $\Gamma(p,t)$ given in Appendix~A.9
of\cite{dpp2}, namely
  \begin{eqnarray*}
&\Gamma(p,t)=\left\langle\,R_1, R_2, R_3, J\vphantom{(R_1R_2R_3R_2^{-1})^{\frac{6p}{p-6}}}\,\right.\left\vert\,  
     R_1^p,\, J^3,\, (R_1J)^{2k},\, JR_1J^{-1}=R_2,\, JR_2J^{-1}=R_3,\,  
     \right. & \\ 
& \left. \br_3(R_1,R_2),\,(R_1R_2)^{\frac{6p}{p-6}}, 
     (JR_2R_1)^{\frac{4kp}{(2k-4)p-4k}}
     \,\right\rangle&
   \end{eqnarray*}
where the last two relations are omitted when the corresponding
exponent is negative.

We wish to construct a homomorphism of $\Gamma(p,t)$ onto the
$(2,3,p)$ triangle group
$$
T_{2,3,p}=\FPG{\beta,\gamma}{(\beta\gamma)^2,\beta^3,\gamma^p}.
$$
For convenience, recall that the optimal index of a torsion-free
subgroup in the hyperbolic triangle groups are as in
Table~\ref{tab:index_23p} (cf.~\cite{eek}).
\begin{table}
\begin{tabular}{c|cccccc}
  $p$       & 7  & 8  & 9  & 10 & 12 & 18 \\
  \hline
  Index     & 84 & 48 & 36 & 30 & 24 & 18
\end{tabular}
\caption{Optimal index of torsion-free subgroups in the triangle group
  $T_{2,3,p}$.}\label{tab:index_23p}
\end{table}

\begin{prop}\label{prop:basic-23p}
  Suppose $p>6$ and $\Gamma(p,t)$ satisfies the Mostow $\Sigma$-INT
  condition. Then there is a unique group homomorphism
  $\varphi:\Gamma(p,t)\rightarrow T_{2,3,p}$ such that $\varphi(J)=\beta$
  and $\varphi(R_1)=\gamma$.
\end{prop}

\begin{pf}
  If the map exists, then we must have
  $\varphi(R_2)=\beta\gamma\beta^{-1}$ and
  $\varphi(R_3)=\beta^{-1}\gamma\beta$. We then check that, under the
  hypotheses of Proposition~\ref{prop:basic-23p}, for every relator
  $w_1\dots w_n$ in the presentation of $\Gamma(p,t)$, the relation
  $\varphi(w_1)\dots\varphi(w_n)=id$ holds in $T_{2,3,p}$ (in that case, we
  say that the relation is compatible).

  For example,
  $\varphi(R_1)\varphi(R_2)\varphi(R_1)=\gamma\beta\gamma\beta^{-1}\gamma=\beta^{-2}\gamma=\beta\gamma$,
  and
  $\varphi(R_2)\varphi(R_1)\varphi(R_2)=\beta\gamma\beta^{-1}\gamma\beta\gamma\beta^{-1}=\beta\gamma\beta^2\gamma\beta\gamma\beta^2=\gamma^{-1}\beta^{-1}=\beta\gamma$,
  so the braid relation $R_1R_2R_1(R_2R_1R_2)^{-1}$ is compatible.

  For the other relations, one checks that
  $\varphi(R_1)\varphi(R_2)=\beta$ and
  $\varphi(R_1)\varphi(R_2)\varphi(J)=\gamma^{-2}$. It is then easy to
  verify that for the relevant values of $k=2,3,4$ and $5$,
  $6p/(p-6)$ is a multiple of 3 and $4kp/((k-2)p-2k)$ is a multiple of $p$.
\end{pf}

Proposition~\ref{prop:basic-23p} implies
Proposition~\ref{prop:maps-to-curves} for the Mostow groups
$\Gamma(7,\frac{9}{14})$, $\Gamma(8,\frac{5}{8})$,
$\Gamma(9,\frac{11}{18})$, $\Gamma(10,\frac{3}{5})$,
$\Gamma(12,\frac{7}{12})$, $\Gamma(18,\frac{5}{9})$,
$\Gamma(7,\frac{13}{42})$, $\Gamma(8,\frac{7}{24})$,
$\Gamma(9,\frac{5}{18})$, $\Gamma(10,\frac{4}{15})$,
$\Gamma(12,\frac{1}{4})$, $\Gamma(18,\frac{12}{9})$,
$\Gamma(8,\frac{1}{8})$, $\Gamma(12,\frac{1}{12})$, $\Gamma(10,0)$.
Note that the first 6 groups in this list are Livne lattices; one can
check that the corresponding maps are the same as the ones that come
from the Livn\'e construction (the construction is conveniently
reviewed in~\cite{deligne-mostow-book}, \S 16).

For the other cases of Proposition~\ref{prop:maps-to-curves}, we use
$\Gamma_\mu$ instead of $\Gamma_{\mu,\Sigma}$. We define two kinds of
subgroups of $\Gamma_\mu$, namely
\begin{equation}
  \label{eq:subgroups}
  K = \langle\langle A_1,A_2,A_3\rangle\rangle,\quad L = \langle\langle A_1,B_2,B_3\rangle\rangle,
\end{equation}
see the notation in equation~\eqref{eq:generators-Gamma_mu}.

We then use group theory software to compute a simplified presentation
for $\Gamma_\mu/K$ and $\Gamma_\mu/L$, and find that in each of the
cases listed in Proposition~\ref{prop:maps-to-curves}, the
corresponding quotient is a $p,q,r$ triangle group.

The details are listed in Table~\ref{tab:triangle-groups}. In the last
column of the table, we list the index of the subgroups we obtain that
map to a non-Abelian free group. The factor 6 comes from the index of
$\Gamma_\mu$ in $\Gamma_{\mu,\Sigma}$, the other factor comes
from the minimal index of a torsion-free subgroup in a
Fuchsian group~\cite{eek}.
\begin{table}[htbp]
  \centering
  \begin{tabular}[htbp]{|r|c|c|c|c|}
    \hline
    $\Gamma(p,t)$     &  Cocompact, Arithmetic?  & $\Gamma_\mu/K$ & $\Gamma_\mu/L$  & Index \\
    \hline
    $\Gamma(8,5/8)$   &     C,A                  & 4,4,4           & 2,2,4          &  $6\cdot 8 = 48$\\
    $\Gamma(10,3/5)$  &     C,A                  & 5,5,5           & 2,2,5          &  $6\cdot 5 = 30$\\
    $\Gamma(12,7/12)$ &     C,A                  & 6,6,6           & 2,2,2          &  $6\cdot 12 = 72$\\
    $\Gamma(18,5/9)$  &     C,A                  & 9,9,9           & 2,2,3          &  $6\cdot 9 = 54$\\
    $\Gamma(8,7/24)$  &     C,NA                 & 4,4,4           & 3,3,4          &  $6\cdot 8 = 48$\\
    $\Gamma(10,4/15)$ &     C,NA                 & 5,5,5           & 3,3,5          &  $6\cdot 5 = 30$\\
    $\Gamma(12,1/4)$  &     C,A                  & 6,6,6           & 2,3,3          &  $6\cdot 12 = 72$\\
    $\Gamma(18,2/9)$  &     C,A                  & 9,9,9           & 3,3,3          &  $6\cdot 3 = 18$\\
    $\Gamma(6,1/6)$   &     NC,NA                & 3,3,3           & 3,4,4          &  $6\cdot 3 = 18$\\
    $\Gamma(8,1/8)$   &     C,A                  & 4,4,4           & 4,4,4          &  $6\cdot 12 = 72$\\
    $\Gamma(12,1/12)$ &     C,A                  & 6,6,6           & 2,2,2          &  $6\cdot 12 = 72$\\
    $\Gamma(4,3/20)$  &     C,A                  & 2,2,2           & 2,5,5          &  $6\cdot 20 = 120$\\
    $\Gamma(10,0)$    &     C,A                  & 5,5,5           & 5,5,5          &  $6\cdot 5 = 30$\\
    $\Gamma(4,1/12)$  &     C,NA                 & 2,2,2           & 2,6,6          &  $6\cdot 12 = 72$\\
    $\Gamma(6,0)$     &     NC,A                 & 3,3,3           & 3,6,6          &  $6\cdot 3 = 18$\\
    $\Gamma(4,0)$     &     C,A                  & 2,2,2           & 2,8,8          &  $6\cdot 8 = 48$\\
    \hline
  \end{tabular}
  \caption{List of maps to triangle groups; in the last column, we
    list the index of the subgroups we find that map to a non-Abelian
    free group.}\label{tab:triangle-groups}
\end{table}

Recall that when $1/p+1/q+1/r$, the $(p,q,r)$ is a lattice in
$PSL(2,\R)$, and this is satisfied for at least one triple $(p,q,r)$
of each row of table~\ref{tab:triangle-groups}, so the
proof of Proposition~\ref{prop:maps-to-curves} is complete.

Now Lemma~\ref{lem:sl2r} implies the following.
\begin{thm}\label{thm:large}
  Let $\Gamma$ be one of the Mostow lattices listed in
  Proposition~\ref{prop:maps-to-curves}. Then $\Gamma$ is large.
\end{thm}
Note also that our proof give a way to describe explicit subgroups,
and explicit morphisms to $F_2$.

\subsection{Neat subgroups}\label{sec:neat}

We summarize the results of our search for neat subgroups of small
index in Tables~\ref{tab:neat-3s1} through~\ref{tab:neat-m3-0}. In
these tables, for each triangle group, we give
\begin{itemize}
\item The Euler characteristic and the least common multiple (LCM) of
  the orders of isotropy groups; this gives a lower bound for the
  index of torsion-free subgroups. Note that by Noether's formula
  (see section~\ref{sec:lattices}), compact ball quotients have Euler
  characteristic an integer multiple of $3$. In some cases, this
  yields a slightly higher lower bound (given by $3\cdot$LCM);
\item The arithmeticity (A/NA) and
  cocompactness (C/NC) of the group;
\item Its Abelianization;
\item Generators for (representatives of) its cusps (inside one given
  box corresponding to the description of the cusps, each line
  corresponds to different conjugacy classes of cusps).
\item The order of the smallest congruence image such that the
  corresponding principal congruence subgroup is torsion-free (we also
  an explicit description of the congruence image group, if we know
  one). In some cases, the congruence image seems too large to compute
  anything about it, in which case we simply write ``?'' in the
  corresponding box of the table. In the last column, we list the
  rational prime whose prime ideal factor was used to get a
  torsion-free congruence subgroup (in some rare cases, one can get a
  smaller congruence by reduction modulo a non-prime ideal). 
\item If we found any neat subgroup with index smaller than the order
  of the congruence image, we list some basic invariants for the
  corresponding subgroups (index, index of the normal core, quotient
  by the normal core, abelianization, self-intersections of the
  elliptic curves in the toroidal compactification, first Betti
  number).
\end{itemize}

\begin{rmk}
  \begin{enumerate}
  \item When the Congruence image column of the following tables
    contains a question mark ``?'', we mean that the corresponding
    finite group seems too large to compute anything. We suspect that
    one could probably work a little more and compute the order of
    these groups and identify them as explicit projectivized linear
    groups over finite fields, but it is probably hopeless to try and
    compute anything about the corresponding principal congruence
    subgroups (presentation, Abelianization).
  \item In Tables~\ref{tab:neat-3s1} through~\ref{tab:neat-m3-0}, we
    write $\Z_n$ for $\Z/n\Z$, and $PU(m,n)$ denotes $PU(m,\F_n)$,
    $PGL(m,n)$ denotes $PGL(m,\F_n)$ where $F_n$ is the field with $n$
    elements ($n$ is a power of a single prime).
  \end{enumerate}
\end{rmk}

\begin{table}[htbp]
  \centering
  {\large $\Sc(3,\sigma_1)$}\\
\begin{tabular}[htbp]{|c|c|c|c|c|c|c|c|}
    \hline
    $\chi$ & LCM & A? & C? & Ab     &  Cusps    & Cusp relations   & Congr. image (mod 5)\\
    \hline
    $\frac{2}{9}$  &  72 &  NA  &  NC  & $\Z_3$ & $R_1,R_2$ & $\br_6(R_1,R_2)$   &  378\,000, $PGU(3,5)$\\
    \hline
  \end{tabular}
  \begin{tabular}[htbp]{|c|c|c|c|c|c|c|}
    \hline
    Index           &  Core Index & Quotient                & Ab                          & Self-int           & $b_1$ \\
    \hline
    $432=6\cdot 72$ &   16868     & $\Z_3\times PGL(3,3)$ & $\Z_3 \oplus \Z_6 \oplus \Z_{18}$ & $(-3)^6, (-9)^{14}$ & $0$\\
    \hline
    $864=12\cdot 72$ &   18144     & $\Z_3\times PSU(3,3)$ & $\Z_2^2\oplus\Z^2$              & $(-2)^{36}$          & $2$\\
    \hline
  \end{tabular}
  \caption{Neat subgroups of $\Sc(3,\sigma_1)$ with core index $\leq 20000$}
  \label{tab:neat-3s1}
\end{table}

\begin{table}[htbp]
  \centering
  {\large $\Sc(4,\sigma_1)$}\\
  \begin{tabular}[htbp]{|c|c|c|c|c|c|c|c|}
    \hline
   $\chi$& LCM & A? & C? & Ab        &  Cusps               & Cusp relations               & Congr. image (mod 3)\\
    \hline
    $\frac{7}{16}$  &  96 &  NA  &  NC  & $\Z_4$ & $1,23\bar2$ & $\br_4(1,23\bar2)$  & 42\,456\,960, PGL(3,9)\\
    \hline
  \end{tabular}
  \begin{tabular}[htbp]{|c|c|c|c|c|c|c|}
    \hline
    Index          &  Core Index & Quotient                & Ab                                   & Self-int            & $b_1$ \\
    \hline
    96             &   96        &                         & $\Z_2^2 \oplus \Z_4^3 \oplus \Z^4$        & $(-8)^{24}$         & $4$\\
    \hline
    96             &   384       &                         & $\Z_2^5 \oplus \Z_8^2$               & $(-8)^{6}$          & $0$\\
    \hline
    96             &   384       &                         & $\Z_2^5 \oplus \Z_4^3$               & $(-4)^{12}$         & $0$\\
    \hline
    96             &   384       &                         & $\Z_2^2 \oplus \Z_4^3 \oplus \Z_8^2$ & $(-4)^{12}$         & $0$\\
    \hline
    96             &   384       &                         & $\Z_2^5 \oplus \Z_4^2 \oplus \Z_8^2$ & $(-2)^{8},(-4)^8$ & $0$\\
    \hline
    96             &   384       &                         & $\Z_2^3 \oplus \Z_4^3 \oplus \Z^2$   & $(-2)^{8},(-4)^8$ & $2$\\
    \hline
  \end{tabular}
  \caption{Neat subgroups of $\Sc(3,\sigma_1)$ with core index $\leq 384$}
  \label{tab:neat-4s1}
\end{table}

\begin{table}[htbp]
  \centering
  {\large $\Sc(6,\sigma_1)$}\\
  \begin{tabular}[htbp]{|c|c|c|c|c|c|c|c|}
    \hline
   $\chi$ & LCM & A? & C? & Ab        &  Cusps               & Cusp relations                                    & Congr. image (mod 5)\\
    \hline
    $\frac{43}{72}$ &  72 &  NA  &  NC  & $\Z_6$ & $1,232\bar3\bar2$ & $\br_3(1,232\bar3\bar2)$  & 378\,000, $PGU(3,5)$\\
          &     &      &      &          & $1,\bar3\bar2323$ & $\br_3(1,\bar3\bar2323)$  &   \\
    \hline
  \end{tabular}
  \begin{tabular}[htbp]{|c|c|c|c|c|c|}
    \hline
    Index               &  Core Index & Quotient                & Ab                                  & Self-int            & $b_1$ \\
    \hline
    $1728=24\cdot 72$   & 36288       & $PSU(3,3)\times\Z_3$  & $\Z_2^3\oplus \Z_4\oplus \Z^4$  & $(-2)^{72}$         & $4$\\
    \hline
  \end{tabular}
  \caption{Neat subgroups of $\Sc(6,\sigma_1)$.}
  \label{tab:neat-6s1}
\end{table}

\begin{table}[htbp]
  \centering
  {\large $\Sc(3,\bar \sigma_4)$}\\
  \begin{tabular}[htbp]{|c|c|c|c|c|c|c|c|}
    \hline
   $\chi$ & LCM & A? & C? & Ab        &  Cusps               & Cusp relations &  Congr. image (mod 5)\\
    \hline
    $\frac{2}{63}$  & 504 &  A   &  C   & $\Z_3$ &                       &                &   378\,000, $PGU(3,5)$\\
    \hline
  \end{tabular}
  \begin{tabular}[htbp]{|c|c|c|c|c|c|}
    \hline
    Index               &  Core Index & Quotient                & Ab                                  & Self-int            & $b_1$ \\
    \hline
    $6048=12\cdot 504$  &  18144       & $PSU(3,3)\times\Z_3$  & $\Z^{10}$                          &                     & $10$\\
    \hline
  \end{tabular}
  \caption{Neat subgroups of $\Sc(3,\bar\sigma_4)$.}
  \label{tab:neat-3s4c}
\end{table}

\begin{table}[htbp]
  \centering
  {\large $\Sc(4,\bar \sigma_4)$}\\
  \begin{tabular}[htbp]{|c|c|c|c|c|c|c|c|}
    \hline
   $\chi$ & LCM & A? & C? & Ab         &  Cusps               & Cusp relations    & Congr. image (mod 3)\\
    \hline
    $\frac{25}{224}$  & 672 &  NA  &  NC  & $\{0\}$  & $R_1,R_2$             & $\br_4(R_1,R_2)$ &  6\,048, $PGU(3,3)$\\
    \hline
  \end{tabular}
  \begin{tabular}[htbp]{|c|c|c|c|c|c|}
    \hline
    Index               &  Core Index & Quotient                & Ab                                  & Self-int            & $b_1$ \\
    \hline
    $6048=9\cdot 672$   &  6048       & $PSU(3,3)$              & $\Z_3^{8}$                          & $(-6)^{56}$           & $0$\\
    \hline
  \end{tabular}
  \caption{Neat subgroups of $\Sc(3,\bar\sigma_4)$.}
  \label{tab:neat-4s4c}
\end{table}

\begin{table}[htbp]
  \centering
  {\large $\Sc(5,\bar \sigma_4)$}\\
  \begin{tabular}[htbp]{|c|c|c|c|c|c|c|c|}
    \hline
   $\chi$   & LCM  & A? & C? & Ab         &  Cusps               & Cusp relations  & Congr. image (mod 3)\\
    \hline
    $\frac{47}{280}$  & 4200 &  NA  &  C  & $\{0\}$     &                      &                 &   ?\\
    \hline
  \end{tabular}
  \caption{Invariants for $\Sc(5,\bar\sigma_4)$}
  \label{tab:5s4c}
\end{table}

\begin{table}[htbp]
  \centering
  {\large $\Sc(6,\bar \sigma_4)$}\\
  \begin{tabular}[htbp]{|c|c|c|c|c|c|c|c|}
    \hline
   $\chi$   & LCM  & A? & C? & Ab         &  Cusps                 & Cusp relations               & Congr. image (mod 5)\\
    \hline
    $\frac{25}{126}$  & 252  &  NA  &  NC  & $\Z_3$     & $1,23\bar2$   & $\br_3(1,23\bar2)$  & 378\,000, $PGU(3,5)$\\
    \hline
  \end{tabular}
  \begin{tabular}[htbp]{|c|c|c|c|c|c|}
    \hline
    Index               &  Core Index & Quotient                & Ab                                  & Self-int            & $b_1$ \\
    \hline
   $3024=12\cdot 252$   & $378\,000$     & $PGU(3,5)$           & $\Z_5^3$                            &  $(-1)^4,(-5)^4$    & 0\\
    \hline
  \end{tabular}
  \caption{Torsion-free subgroups of $\Sc(6,\bar\sigma_4)$.}
  \label{tab:neat-6s4c}
\end{table}

\begin{table}[htbp]
  \centering
  {\large $\Sc(8,\bar \sigma_4)$}\\
  \begin{tabular}[htbp]{|c|c|c|c|c|c|c|c|}
    \hline
   $\chi$   & LCM  & A? & C? & Ab         &  Cusps                 & Cusp relations   & Congr. image (mod 3)\\
    \hline
    $\frac{99}{448}$  & 1344 &  NA  &  C   & $\{0\}$    &                        &                  & 42\,456\,960, PGL(3,9)\\
    \hline
  \end{tabular}
  \begin{tabular}[htbp]{|c|c|c|c|c|c|}
    \hline
    Index                    & Core Index & Quotient                & Ab                                  & Self-int            & $b_1$ \\
    \hline
   $131712 = 98\cdot 1344$   & 5663616  & $PGU(3,7)$              & ?                                   &                     & ?\\
    \hline
  \end{tabular}
  \caption{Torsion-free subgroups of $\Sc(8,\bar\sigma_4)$.}
  \label{tab:neat-8s4c}
\end{table}

\begin{table}[htbp]
  \centering
  {\large $\Sc(12,\bar \sigma_4)$}\\
  \begin{tabular}[htbp]{|c|c|c|c|c|c|c|c|}
    \hline
   $\chi$     & LCM  & A? & C? & Ab         &  Cusps                 & Cusp relations  & Congr. image (mod 5)\\
    \hline
    $\frac{221}{1008}$  & 1008 &  NA  &  C   & $\Z_3$   &                        &                 &   ? \\
    \hline
  \end{tabular}
  \begin{tabular}[htbp]{|c|c|c|c|c|c|}
    \hline
    Index                    & Core Index & Quotient                & Ab                                  & Self-int            & $b_1$ \\
    \hline
   $6048 = 6\cdot 1008$      & 18144      & $PGU(3,3)\times \Z_3$   & $\Z_3^4$                            &                     & $0$\\
    \hline
  \end{tabular}
  \caption{Torsion-free subgroups of $\Sc(12,\bar\sigma_4)$ with core
    index $\leq 20000$.}
  \label{tab:neat-12s4c}
\end{table}

\begin{table}[htbp]
  \centering
  {\large $\Sc(2,\sigma_5)$}\\
  \begin{tabular}[htbp]{|c|c|c|c|c|c|c|c|}
    \hline
   $\chi$     & LCM  & A? & C? & Ab         &  Cusps                 & Cusp relations  & Congr. image (mod 5)\\
    \hline
    $\frac{1}{45}$      & 360  &  A   &  C   & $\Z_6$   &                        &                 &      378\,000, $PGU(3,5)$\\
    \hline
  \end{tabular}
  \begin{tabular}[htbp]{|c|c|c|c|c|c|}
    \hline
    Index                    & Core Index & Quotient                & Ab                                    & Self-int            & $b_1$ \\
    \hline
   $54,000=50\cdot 1080$      & 378\,000  & $PGU(3,5)$              & $\Z_2^{10}\oplus\Z_{10}\oplus\Z_{70}$ &                     & $0$\\
    \hline
  \end{tabular}
  \caption{Torsion-free subgroups of $\Sc(2,\sigma_5)$}
  \label{tab:neat-2s5}
\end{table}

\begin{table}[htbp]
  \centering
  {\large $\Sc(3,\sigma_5)$}\\
  \begin{tabular}[htbp]{|c|c|c|c|c|c|c|c|}
    \hline
   $\chi$     & LCM  & A? & C? & Ab         &  Cusps                                                              & Cusp relations                & Congr. image (mod 5)\\
    \hline 
    $\frac{49}{180}$    & 360  &  NA  &  NC  & $\Z_3^2$   & {\scriptsize $c_1=2, c_2=2(\bar J\bar1)^2121\bar2\bar1(\bar J\bar1)^32$} & {\scriptsize $c_1^3,c_2^6,\br_4(c_1,c_2)$} &      378\,000, $PGU(3,5)$\\
    \hline
  \end{tabular}
  \begin{tabular}[htbp]{|c|c|c|c|c|c|}
    \hline
    Index                    & Core Index & Quotient                & Ab                                    & Self-int            & $b_1$ \\
    \hline
   $360$                     & 360        &                         & $\Z_3\oplus\Z^8$                      & $(-1)^{20}$         & $8$\\
    \hline
   $360$                     & 360        &                         & $\Z_3^9$                              & $(-1)^{20}$         & $0$\\
    \hline
   $360$                     & 1080       &                         & $\Z_3^2\oplus\Z^6$                    & $(-1)^{20}$         & $6$\\
    \hline
  \end{tabular}
  \caption{Neat subgroups of $\Sc(3,\sigma_5)$ with core index $\leq 10,000$}
  \label{tab:neat-3s5}
\end{table}

\begin{table}[htbp]
  \centering
  {\large $\Sc(4,\sigma_5)$}\\
  \begin{tabular}[htbp]{|c|c|c|c|c|c|c|c|}
    \hline
   $\chi$     & LCM  & A? & C? & Ab         &  Cusps    & Cusp relations   & Congr. image (mod 5)\\
    \hline 
    $\frac{17}{36}$     & 720  &  NA  &  NC  & $\Z_6$   & $R_1,R_2$ & $\br_4(R_1,R_2)$ & 152\,334\,000\,000, $PGL(3,25)$\\
    \hline
  \end{tabular}
  \caption{Invariants for $\Sc(4,\sigma_5)$}
  \label{tab:4s5}
\end{table}

\begin{table}[htbp]
  \centering
  {\large $\Sc(3,\sigma_{10})$}\\
  \begin{tabular}[htbp]{|c|c|c|c|c|c|c|c|}
    \hline
   $\chi$     & LCM  & A? & C? & Ab         &  Cusps    & Cusp relations   & Congr. image (mod 5)\\
    \hline 
    $\frac{1}{45}$      & 360  &  A   &   C  & $\Z_3$   &           &                  & 378\,000, $PGU(3,5)$\\
    \hline
  \end{tabular}
  \begin{tabular}[htbp]{|c|c|c|c|c|c|}
    \hline
    Index                    & Core Index & Quotient                & Ab                                    & Self-int            & $b_1$ \\
    \hline
     $2160=2\cdot 1080$      &  2160      &                         & $\Z^{16}$                             &                     & $16$\\
    \hline
     $2160=2\cdot 1080$      &  34560     &                         & $\Z_2^8\oplus\Z_8$                    &                     & $0$\\
    \hline
     $2160=2\cdot 1080$      &  34560     &                         & $\Z_2^2\oplus\Z_4^2\oplus\Z_{16}$      &                    & $0$\\
    \hline
  \end{tabular}
  \caption{Neat subgroups of $\Sc(3,\sigma_{10})$ with core index $\leq 40,000$}
  \label{tab:neat-3s10}
\end{table}

\begin{table}[htbp]
  \centering
  {\large $\Sc(4,\sigma_{10})$}\\
  \begin{tabular}[htbp]{|c|c|c|c|c|c|c|c|}
    \hline
   $\chi$     & LCM  & A? & C? & Ab         &  Cusps    & Cusp relations   & Congr. image (mod 3)\\
    \hline 
    $\frac{3}{32}$      & 480  &  A   &   C  & $\{1\}$    &           &                  & 42\,456\,960, $PGL(3,9)$\\
    \hline
  \end{tabular}
  \begin{tabular}[htbp]{|c|c|c|c|c|c|}
    \hline
    Index                    & Core Index & Quotient                & Ab                                    & Self-int            & $b_1$ \\
    \hline
     $12,000=25\cdot 480$    &  $372,000$ & $PGL(3,5)$              & $\Z_5\oplus\Z_{155}$                  &                     & $0$\\
    \hline
  \end{tabular}
  \caption{Neat subgroups of $\Sc(4,\sigma_{10})$}
  \label{tab:neat-4s10}
\end{table}

\begin{table}[htbp]
  \centering
  {\large $\Sc(5,\sigma_{10})$}\\
  \begin{tabular}[htbp]{|c|c|c|c|c|c|c|c|}
    \hline
   $\chi$     & LCM  & A? & C? & Ab         &  Cusps    & Cusp relations   & Congr. image (mod 3)\\
    \hline 
    $\frac{1}{8}$       & 600  &  A   &   C  & $\Z_5$   &           &                  & 42\,573\,600, $PGU(3,9)$\\
    \hline
  \end{tabular}
  \begin{tabular}[htbp]{|c|c|c|c|c|c|}
    \hline
    Index                    & Core Index & Quotient                & Ab                                    & Self-int            & $b_1$ \\
    \hline
     $600$                   &  600       &                         & $\Z^{20}$                             &                     & $20$\\
    \hline
     $600$                   &  15\,000    &                         & $\Z_5^{7}$                            &                     & $0$\\
    \hline
     $600$                   &  15\,000    &                         & $\Z_5\oplus\Z^{8}$                    &                     & $8$\\
    \hline
  \end{tabular}
  \caption{Neat subgroups of $\Sc(5,\sigma_{10})$ with core index $\leq 20,000$}
  \label{tab:neat-5s10}
\end{table}

\begin{table}[htbp]
  \centering
  {\large $\Sc(10,\sigma_{10})$}\\
  \begin{tabular}[htbp]{|c|c|c|c|c|c|c|c|}
    \hline
   $\chi$     & LCM  & A? & C? & Ab         &  Cusps    & Cusp relations   & Congr. image (mod 3)\\
    \hline 
    $\frac{13}{100}$    & 300  &  A   &   C  & $\Z_5$   &           &                  & 42\,573\,600, $PGU(3,9)$\\
    \hline
  \end{tabular}
  \begin{tabular}[htbp]{|c|c|c|c|c|c|}
    \hline
    Index                    & Core Index & Quotient                & Ab                                    & Self-int            & $b_1$ \\
    \hline
     $300$                   &  18,000    &                         & $\Z_5^2\oplus\Z_{15}$                 &                     & $0$\\
    \hline
  \end{tabular}
  \caption{Neat subgroups of $\Sc(10,\sigma_{10})$ with core index $\leq 20,000$}
  \label{tab:neat-10s10}
\end{table}

\begin{table}[htbp]
  \centering
  {\large $\Tc(3,{\bf S}_2)$}\\
  \begin{tabular}[htbp]{|c|c|c|c|c|c|c|c|}
    \hline
   $\chi$     & LCM  & A? & C? & Ab         &  Cusps    & Cusp relations   & Congr. image (mod 5)\\
    \hline 
    $\frac{2}{15}$      & 360  &  A   &   C  & $\Z_3$   &           &                  & 378\,000, $PGU(3,5)$\\
    \hline
  \end{tabular}
  \begin{tabular}[htbp]{|c|c|c|c|c|c|}
    \hline
    Index                    & Core Index & Quotient                & Ab                                    & Self-int            & $b_1$ \\
    \hline
     $360$                   &  360       &                         & $\Z^{16}$                             &                     & $16$\\
    \hline
  \end{tabular}
  \caption{Neat subgroups of $\Tc(3,{\bf S}_2)$ with core index $\leq 2,000$}
  \label{tab:neat-3S2}
\end{table}

\begin{table}[htbp]
  \centering
  {\large $\Tc(4,{\bf S}_2)$}\\
  \begin{tabular}[htbp]{|c|c|c|c|c|c|c|c|}
    \hline
   $\chi$     & LCM  & A? & C? & Ab         &  Cusps    & Cusp relations   & Congr. image (mod 5)\\
    \hline 
    $\frac{1}{3}$       & 480  &  NA  &   NC & $\{1\}$    & $R_1,R_2$ & $\br_4(R_1,R_2)$ & 42\,456\,960, $PGL(3,9)$\\
    \hline
  \end{tabular}
  \begin{tabular}[htbp]{|c|c|c|c|c|c|}
    \hline
    Index                    & Core Index & Quotient                & Ab                                    & Self-int            & $b_1$ \\
    \hline
     $155,520=324\cdot 480$  & 42456960   &                         & $?$                                   &                     & $?$\\
    \hline
  \end{tabular}
  \caption{Neat subgroups of $\Tc(4,{\bf S}_2)$}
  \label{tab:neat-4S2}
\end{table}

\begin{table}[htbp]
  \centering
  {\large $\Tc(5,{\bf S}_2)$}\\
  \begin{tabular}[htbp]{|c|c|c|c|c|c|c|c|}
    \hline
   $\chi$     & LCM  & A? & C? & Ab         &  Cusps    & Cusp relations   & Congr. image (mod 3)\\
    \hline 
    $\frac{133}{300}$   & 600  &  NA  &   C  & $\Z_5$   &           &                  & 23\,619\,600\\
    \hline
  \end{tabular}
  \begin{tabular}[htbp]{|c|c|c|c|c|c|}
    \hline
    Index                    & Core Index & Quotient                & Ab                                    & Self-int            & $b_1$ \\
    \hline
     $3600=2\cdot 1800$      &    3600    &                         & $\Z_3^8$                              &                     & $0$\\
    \hline
  \end{tabular}
  \caption{Neat subgroups of $\Tc(5,{\bf S}_2)$}
  \label{tab:neat-5S2}
\end{table}

\begin{table}[htbp]
  \centering
  {\large $\Tc(3,{\bf E}_2)$}\\
  \begin{tabular}[htbp]{|c|c|c|c|c|c|c|c|}
    \hline
   $\chi$     & LCM  & A? & C? & Ab         &  Cusps    & Cusp relations   & Congr. image (mod 5)\\
    \hline 
    $\frac{1}{4}$       & 72   &  A  &   NC  & $\Z_3^2$   & {\scriptsize $c_1=2,c_2=(\bar3\bar2\bar1)^2\bar2\bar1 2$} & {\scriptsize $c_1^3, c_2^6(c_1c_2)^{-4}, (c_1c_2)^2\textrm{ central}$} & 378\,000, $PGU(3,5)$\\
    \hline
  \end{tabular}
  \begin{tabular}[htbp]{|c|c|c|c|c|c|}
    \hline
    Index                    & Core Index & Quotient                & Ab                                    & Self-int            & $b_1$ \\
    \hline
     $72$                    &    72      &                         & $\Z^8$                                & $(-1)^{12}$         & $8$\\
    \hline
     $72$                    &    648     &                         & $\Z_3^2\oplus\Z^4$                    & $(-1)^3,(-3)^3$     & $4$\\
    \hline
     $72$                    &    648     &                         & $\Z_3^3\oplus\Z^2$                    & $(-3)^4$            & $2$\\
    \hline
     $72$                    &    648     &                         & $\Z_3^4\oplus\Z^2$                    & $(-1)^3,(-3)^3$     & $2$\\
    \hline
  \end{tabular}
  \caption{Neat subgroups of $\Tc(3,{\bf E}_2)$ of core index $\leq 1000$}
  \label{tab:neat-3E2}
\end{table}

\begin{table}[htbp]
  \centering
  {\large $\Tc(4,{\bf E}_2)$}\\
  \begin{tabular}[htbp]{|c|c|c|c|c|c|c|c|}
    \hline
   $\chi$     & LCM  & A? & C? & Ab               &  Cusps    & Cusp relations   & Congr. image (mod 5)\\
    \hline 
    $\frac{17}{32}$     & 96   &  NA &   NC  & $\Z_2\oplus\Z_4$ & $R_1,R_2$ & $\br_4(R_1,R_2)$ & 50\,778\,000\,000, $PSL(3,25)$\\
              &      &     &       &                  & $R_1,R_3$ & $\br_4(R_1,R_3)$ & \\
              &      &     &       &                  & $R_1,R_2R_3R_2^{-1}$ & $\br_4(R_1,R_2R_3R_2^{-1})$ & \\
    \hline
  \end{tabular}
  \begin{tabular}[htbp]{|c|c|c|c|c|c|}
    \hline
    Index                    & Core Index & Quotient                & Ab                                    & Self-int            & $b_1$ \\
    \hline
     $96$                    &    192     &                         & $\Z_2^4\oplus\Z_4^2\oplus\Z^4$        & $(-2)^{18}$         & $4$\\
    \hline
     $96$                    &    768     &                         & $\Z_2^4\oplus\Z_4^2\oplus\Z^4$        & $(-2)^{10},(-4)^{4}$         & $4$\\
    \hline
     $96$                    &    768     &                         & $\Z_2^6\oplus\Z_4^3\oplus\Z_8$        & $(-2)^{10},(-4)^{4}$         & $0$\\
    \hline
  \end{tabular}
  \caption{Neat subgroups of $\Tc(4,{\bf E}_2)$ of core index $\leq 1000$}
  \label{tab:neat-4E2}
\end{table}

\begin{table}[htbp]
  \centering
  {\large $\Tc(6,{\bf E}_2)$}\\
  \begin{tabular}[htbp]{|c|c|c|c|c|c|c|c|}
    \hline
   $\chi$     & LCM  & A? & C? & Ab               &  Cusps    & Cusp relations   & Congr. image (mod 5)\\
    \hline 
    $\frac{3}{4}$       & 36   &  A   &   NC  & $\Z_6^2$        & $R_2,R_3$ & $\br_3(R_2,R_3)$ & 378\,000, $PGU(3,5)$\\
    \hline
  \end{tabular}
  \begin{tabular}[htbp]{|c|c|c|c|c|c|}
    \hline
    Index                    & Core Index & Quotient                & Ab                                    & Self-int            & $b_1$ \\
    \hline
     $72$                    &     72     &                         & $\Z_2^3\oplus\Z^{8}$                  & $(-1)^{12}$         & $8$\\
    \hline
  \end{tabular}
  \caption{Neat subgroups of $\Tc(6,{\bf E}_2)$ of core index $\leq 200$}
  \label{tab:neat-6E2}
\end{table}

\begin{table}[htbp]
  \centering
  {\large $\Tc(2,{\bf H}_1)$}\\
  \begin{tabular}[htbp]{|c|c|c|c|c|c|c|c|}
    \hline
   $\chi$     & LCM  & A? & C? & Ab               &  Cusps    & Cusp relations   & Congr. image (mod 5)\\
    \hline 
    $\frac{1}{49}$      & 1176 &  A   &   C  & $\Z_2$           &           &                  & ?\\
    \hline
  \end{tabular}
  \caption{Invariants for $\Tc(2,{\bf H}_1)$}
  \label{tab:2H1}
\end{table}

\begin{table}[htbp]
  \centering
  {\large $\Tc(2,{\bf H}_2)$}\\
  \begin{tabular}[htbp]{|c|c|c|c|c|c|c|c|}
    \hline
   $\chi$     & LCM  & A? & C? & Ab               &  Cusps    & Cusp relations   & Congr. image (mod 3)\\
    \hline 
    $\frac{1}{100}$    & 300  &  A   &   C  & $\{1\}$          &           &                  & 42\,573\,600, $PGU(3,9)$\\
    \hline
  \end{tabular}
  \caption{Invariants for $\Tc(2,{\bf H}_2)$}
  \label{tab:2H2}
\end{table}

\begin{table}[htbp]
  \centering
  {\large $\Tc(3,{\bf H}_2)$}\\
  \begin{tabular}[htbp]{|c|c|c|c|c|c|c|c|}
    \hline
   $\chi$     & LCM  & A? & C? & Ab               &  Cusps    & Cusp relations  & Congr. image (mod 7)\\
    \hline 
    $\frac{26}{75}$     & 1800 &  NA  &   C  & $\Z_3$        &           &                  &   ?\\
    \hline
  \end{tabular}
  \caption{Invariants for $\Tc(3,{\bf H}_2)$}
  \label{tab:3H2}
\end{table}

\begin{table}[htbp]
  \centering
  {\large $\Tc(5,{\bf H}_2)$}\\
  \begin{tabular}[htbp]{|c|c|c|c|c|c|c|c|}
    \hline
   $\chi$     & LCM  & A? & C? & Ab               &  Cusps    & Cusp relations   & Congr. image (mod 3)\\
    \hline 
    $\frac{73}{100}$    & 600  &  A   &   C  & $\Z_5$         &           &                  &   23\,619\,600\\
    \hline
  \end{tabular}
  \begin{tabular}[htbp]{|c|c|c|c|c|c|}
    \hline
    Index                    & Core Index & Quotient                & Ab                                    & Self-int            & $b_1$ \\
    \hline
     $600$                   &     600    &                         & $\Z^{20}$                              &                     & $20$\\
    \hline
  \end{tabular}
  \caption{Neat subgroups of $\Tc(5,{\bf H}_2)$}
  \label{tab:neat-5H2}
\end{table}


\begin{table}[htbp]
  \centering
  {\large $\Gamma(5,7/10)$}\\
\begin{tabular}[htbp]{|c|c|c|c|c|c|c|c|}
    \hline
   $\chi$     & LCM  & A? & C? & Ab               &  Cusps    & Cusp relations   & Congr. image (mod 3)\\
    \hline 
    $\frac{1}{200}$     & 600  &  A   &   C  & $\{1\}$          &           &                  &  42\,573\,600, $PGU(3,9)$ \\
    \hline
  \end{tabular}
  \begin{tabular}[htbp]{|c|c|c|c|c|c|}
    \hline
    Index                    & Core Index & Quotient                & Ab                                    & Self-int            & $b_1$ \\
    \hline
     $3000=5\cdot 600$       &  15\,000    &                         & $\Z_5\oplus\Z^4$                      &                     & $4$\\
    \hline
  \end{tabular}
  \caption{Neat subgroups of $\Gamma(5,7/10)$ of core index $\leq 40,000$}
  \label{tab:neat-m5-7.10}
\end{table}

\begin{table}[htbp]
  \centering
  {\large $\Gamma(6,2/3)$}\\
  \begin{tabular}[htbp]{|c|c|c|c|c|c|c|c|}
    \hline
   $\chi$     & LCM  & A? & C? & Ab               &  Cusps    & Cusp relations   & Congr. image (mod 5)\\
    \hline 
    $\frac{1}{72}$      &  72  &  A   &  NC  & $\Z_6$         & $R_1,R_2$ & $\br_3(R_1,R_2)$ &  378\,000, $PGU(3,5)$\\
    \hline
  \end{tabular}
  \begin{tabular}[htbp]{|c|c|c|c|c|c|}
    \hline
    Index                    & Core Index & Quotient                & Ab                                    & Self-int            & $b_1$ \\
    \hline
     $72$                    &    1944    &                         & $\Z_3\oplus\Z^2$                      &  $-1,-3$            & $2$\\
    \hline
     $72$                    &    1944    &                         & $\Z^4$                                &  $(-1)^4$           & $4$\\
    \hline
  \end{tabular}
  \caption{Neat subgroups of $\Gamma(6,2/3)$ of core index $\leq 2000$}
  \label{tab:neat-m6-2.3}
\end{table}

\begin{table}[htbp]
  \centering
  {\large $\Gamma(7,9/14)$}\\
  \begin{tabular}[htbp]{|c|c|c|c|c|c|c|c|}
    \hline
   $\chi$     & LCM  & A? & C? & Ab               &  Cusps    & Cusp relations   & Congr. image (mod 3)\\
    \hline 
     $\frac{1}{49}$     &  588 &  A   &   C  & $\{1\}$          &           &                  &  ? \\
    \hline
  \end{tabular}
  \caption{Invariants for $\Gamma(7,9/14)$}
  \label{tab:m7-9.14}
\end{table}

\begin{table}[htbp]
  \centering
  {\large $\Gamma(8,5/8)$}\\
  \begin{tabular}[htbp]{|c|c|c|c|c|c|c|c|}
    \hline
   $\chi$     & LCM  & A? & C? & Ab               &  Cusps    & Cusp relations   & Congr. image (mod 3)\\
    \hline 
     $\frac{3}{128}$    & 384  &  A   &  C   & $\Z_4$         &           &                  &  42\,456\,960, $PGL(3,9)$ \\
    \hline
  \end{tabular}
  \begin{tabular}[htbp]{|c|c|c|c|c|c|}
    \hline
    Index                    & Core Index & Quotient                & Ab                                    & Self-int            & $b_1$ \\
    \hline
     $466560=1215\cdot 384$  &  42456960  &                         & $?$                                   &                     & $?$\\
    \hline
  \end{tabular}
  \caption{Neat subgroups of $\Gamma(8,5/8)$}
  \label{tab:neat-m8-5.8}
\end{table}

\begin{table}[htbp]
  \centering
  {\large $\Gamma(9,11/18)$}\\
  \begin{tabular}[htbp]{|c|c|c|c|c|c|c|c|}
    \hline
   $\chi$     & LCM  & A? & C? & Ab               &  Cusps    & Cusp relations   & Congr. image (mod 5)\\
    \hline 
    $\frac{2}{81}$      & 324  &   NA &  C   & $\Z_3$         &           &                  &  ? \\
    \hline
  \end{tabular}
  \caption{Invariants for $\Gamma(9,11/18)$}
  \label{tab:m9-11.18}
\end{table}

\begin{table}[htbp]
  \centering
  {\large $\Gamma(10,3/5)$}\\
  \begin{tabular}[htbp]{|c|c|c|c|c|c|c|c|}
    \hline
   $\chi$     & LCM  & A? & C? & Ab               &  Cusps    & Cusp relations   & Congr. image (mod 3)\\
    \hline 
     $\frac{1}{40}$     & 600  &  A   &   C  & $\Z_2$         &           &                  &  42\,573\,600, $PGU(3,9)$ \\
    \hline
  \end{tabular}
  \begin{tabular}[htbp]{|c|c|c|c|c|c|}
    \hline
    Index                    & Core Index & Quotient                & Ab                                    & Self-int            & $b_1$ \\
    \hline
     $583200=972\cdot 600$   & 42573600   &                         & $?$                                    &                     & $?$\\
    \hline
  \end{tabular}
  \caption{Neat subgroups of $\Gamma(10,3/5)$}
  \label{tab:neat-m10-3.5}
\end{table}

\begin{table}[htbp]
  \centering
  {\large $\Gamma(12,7/12)$}\\
  \begin{tabular}[htbp]{|c|c|c|c|c|c|c|c|}
    \hline
   $\chi$     & LCM  & A? & C? & Ab               &  Cusps    & Cusp relations   & Congr. image (mod 5)\\
    \hline 
    $\frac{7}{288}$     & 288  &   A  &  C   & $\Z_{12}$        &           &                  &  152\,334\,000\,000 \\
    \hline
  \end{tabular}
  \caption{Invariants for $\Gamma(12,7/12)$}
  \label{tab:m12-7.12}
\end{table}

\begin{table}[htbp]
  \centering
  {\large $\Gamma(18,5/9)$}\\
  \begin{tabular}[htbp]{|c|c|c|c|c|c|c|c|}
    \hline
   $\chi$     & LCM  & A? & C? & Ab               &  Cusps    & Cusp relations   & Congr. image (mod 5)\\
    \hline 
    $\frac{5}{9}$       & 648  &   A  &   C  & $\Z_6$         &           &                  &  ? \\
    \hline
  \end{tabular}
  \caption{Invariants for $\Gamma(18,5/9)$}
  \label{tab:m18-5.9}
\end{table}

\begin{table}[htbp]
  \centering
  {\large $\Gamma(4,5/12)$}\\
  \begin{tabular}[htbp]{|c|c|c|c|c|c|c|c|}
    \hline
   $\chi$     & LCM  & A? & C? & Ab                 &  Cusps    & Cusp relations   & Congr. image (mod 5)\\
    \hline 
    $\frac{1}{72}$      &  288 &  A   &  C   & $\Z_6$           &           &                  &  152\,334\,000\,000 \\
    \hline
  \end{tabular}
  \begin{tabular}[htbp]{|c|c|c|c|c|c|}
    \hline
    Index                    & Core Index & Quotient                & Ab                                    & Self-int            & $b_1$ \\
    \hline
     $864$                   &   18144    & $PGU(3,3)\times\Z_3$  & $\Z_2^3\oplus\Z^2$                    &                     & $2$\\
    \hline
  \end{tabular}
  \caption{Neat subgroups of $\Gamma(4,5/12)$ of core index $\leq 40,000$}
  \label{tab:neat-m}
\end{table}

\begin{table}[htbp]
  \centering
  {\large $\Gamma(5,11/30)$}\\
  \begin{tabular}[htbp]{|c|c|c|c|c|c|c|c|}
    \hline
   $\chi$     & LCM  & A? & C? & Ab               &  Cusps    & Cusp relations   & Congr. image (mod 7)\\
    \hline 
    $\frac{8}{225}$     & 1800 &   A  &   C  & $\Z_3$         &           &                  &  ? \\
    \hline
  \end{tabular}
  \caption{Invariants for $\Gamma(5,11/30)$}
  \label{tab:m5-11.30}
\end{table}

\begin{table}[htbp]
  \centering
  {\large $\Gamma(6,1/3)$}\\
  \begin{tabular}[htbp]{|c|c|c|c|c|c|c|c|}
    \hline
   $\chi$     & LCM  & A? & C? & Ab               &  Cusps                    & Cusp relations                 & Congr. image (mod 5)\\
    \hline 
    $\frac{1}{18}$      &  18  &   A  &   NC & $\Z_3\oplus\Z_6$ & {\scriptsize $R_1,R_2$}                 & {\scriptsize $\br_3(R_1,R_2)$}               &  378\,000, $PGU(3,5)$\\
              &      &      &      &                  & {\scriptsize $c_1=R_1,c_2=J\bar2\bar1$} & {\scriptsize $c_1^6, c_2^3, \br_4(c_1,c_2)$} &  \\
    \hline
  \end{tabular}
  \begin{tabular}[htbp]{|c|c|c|c|c|c|}
    \hline
    Index                    & Core Index & Quotient                & Ab                                    & Self-int            & $b_1$ \\
    \hline
     $18$                    &    18      &                         & $\Z^4$                                &    $(-1)^4$        & $4$\\
    \hline
     $18$                    &    54      &                         & $\Z_3\oplus\Z^2$                      &    $-1,-3$          & $2$\\
    \hline
     $18$                    &    162     &                         & $\Z^2$                                &    $-1,-3$          & $2$\\
    \hline
  \end{tabular}
  \caption{Neat subgroups of $\Gamma(6,1/3)$ of core index $\leq 1,000$}
  \label{tab:neat-m6-1.3}
\end{table}

\begin{table}[htbp]
  \centering
  {\large $\Gamma(7,13/42)$}\\
  \begin{tabular}[htbp]{|c|c|c|c|c|c|c|c|}
    \hline
   $\chi$     & LCM  & A? & C? & Ab               &  Cusps    & Cusp relations  & Congr. image (mod 5)\\
    \hline 
     $\frac{61}{882}$   & 1764 &  NA  &  C    & $\Z_3$        &           &                 &   ?\\
    \hline
  \end{tabular}
  \caption{Invariants for $\Gamma(7,13/42)$}
  \label{tab:m7-13.42}
\end{table}

\begin{table}[htbp]
  \centering
  {\large $\Gamma(8,7/24)$}\\
  \begin{tabular}[htbp]{|c|c|c|c|c|c|c|c|}
    \hline
   $\chi$     & LCM  & A? & C? & Ab               &  Cusps    & Cusp relations   & Congr. image (mod 5)\\
    \hline 
    $\frac{11}{144}$    & 576  &   NA &   C  & $\Z_6$         &           &                  &  ? \\
    \hline
  \end{tabular}
  \caption{Invariants for $\Gamma(8,7/24)$}
  \label{tab:m8-7.24}
\end{table}

\begin{table}[htbp]
  \centering
  {\large $\Gamma(9,5/18)$}\\
  \begin{tabular}[htbp]{|c|c|c|c|c|c|c|c|}
    \hline
   $\chi$     & LCM  & A? & C? & Ab               &  Cusps    & Cusp relations   & Congr. image (mod 5)\\
    \hline 
     $\frac{13}{162}$   &  324 &   A  &   C  & $\Z_3^2$         &           &                  &  ? \\
    \hline
  \end{tabular}
  \caption{Invariants for $\Gamma(9,5/18)$}
  \label{tab:m9-5.18}
\end{table}

\begin{table}[htbp]
  \centering
  {\large $\Gamma(10,4/15)$}\\
  \begin{tabular}[htbp]{|c|c|c|c|c|c|c|c|}
    \hline
   $\chi$     & LCM  & A? & C? & Ab               &  Cusps    & Cusp relations   & Congr. image (mod 7)\\
    \hline 
    $\frac{37}{450}$    & 450  &  NA  &   C  & $\Z_6$         &           &                  &  ? \\
    \hline
  \end{tabular}
  \caption{Invariants for $\Gamma(10,4/15)$}
  \label{tab:m10-4.15}
\end{table}

\begin{table}[htbp]
  \centering
  {\large $\Gamma(12,1/4)$}\\
  \begin{tabular}[htbp]{|c|c|c|c|c|c|c|c|}
    \hline
   $\chi$     & LCM  & A? & C? & Ab               &  Cusps    & Cusp relations   & Congr. image (mod 5)\\
    \hline 
     $\frac{1}{12}$     & 144  &  A   &  C   & $\Z_3\oplus\Z_6$ &           &                  &  152\,334\,000\,000 \\
    \hline
  \end{tabular}
  \caption{Invariants for $\Gamma(12,1/4)$}
  \label{tab:m12-1.4}
\end{table}

\begin{table}[htbp]
  \centering
  {\large $\Gamma(18,2/9)$}\\
  \begin{tabular}[htbp]{|c|c|c|c|c|c|c|c|}
    \hline
   $\chi$     & LCM  & A? & C? & Ab               &  Cusps    & Cusp relations   & Congr. image (mod 5)\\
    \hline 
    $\frac{13}{162}$    &  162 &  A   &  C   & $\Z_3\oplus\Z_6$ &           &                  &   \\
    \hline
  \end{tabular}
  \caption{Invariants for $\Gamma(18,2/9)$}
  \label{tab:m18-2.9}
\end{table}

\begin{table}[htbp]
  \centering
  {\large $\Gamma(3,1/3)$}\\
  \begin{tabular}[htbp]{|c|c|c|c|c|c|c|c|}
    \hline
   $\chi$     & LCM  & A? & C? & Ab               &  Cusps    & Cusp relations   & Congr. image (mod 5)\\
    \hline 
    $\frac{1}{288}$     &  288 &  A   &  C   & $\Z_3$         &           &                  &  152\,334\,000\,000\\
    \hline
  \end{tabular}
  \begin{tabular}[htbp]{|c|c|c|c|c|c|}
    \hline
    Index                    & Core Index & Quotient                 & Ab                                    & Self-int            & $b_1$ \\
    \hline
     $864$                   &   18,144   & $\Z_3\times PSU(3,3)$    & $\Z^2$                                &                     & $2$\\
    \hline
  \end{tabular}
  \caption{Neat subgroups of $\Gamma(3,1/3)$ of core index $\leq 30,000$}
  \label{tab:neat-m3-1.3}
\end{table}

\begin{table}[htbp]
  \centering
  {\large $\Gamma(4,1/4)$}\\
  \begin{tabular}[htbp]{|c|c|c|c|c|c|c|c|}
    \hline
   $\chi$     & LCM  & A? & C? & Ab               &  Cusps                                   & Cusp relations                              & Congr. image (mod 3)\\
    \hline 
     $\frac{1}{32}$     & 96   &  A   & NC   & $\Z_4$         &  {\scriptsize $c_1=1,c_2=J\bar2\bar1$}  & {\scriptsize $c_1^4, c_2^4, \br_4(c_1,c_2)$} &  6\,048 \\
    \hline
  \end{tabular}
  \begin{tabular}[htbp]{|c|c|c|c|c|c|}
    \hline
    Index                    & Core Index & Quotient                & Ab                                    & Self-int            & $b_1$ \\
    \hline
     $96$                    &   96       &                         & $\Z_2^2\oplus\Z^4$                    & $(-2)^6$            & $4$\\
    \hline
     $96$                    &   384      &                         & $\Z_2\oplus\Z_4^2\oplus\Z^2$          & $(-2)^2,(-4)^2$     & $2$\\
    \hline
  \end{tabular}
  \caption{Neat subgroups of $\Gamma(4,1/4)$ of core index $\leq 1,000$}
  \label{tab:neat-m4-1.4}
\end{table}

\begin{table}[htbp]
  \centering
  {\large $\Gamma(5,1/5)$}\\
  \begin{tabular}[htbp]{|c|c|c|c|c|c|c|c|}
    \hline
   $\chi$     & LCM  & A? & C? & Ab               &  Cusps    & Cusp relations   & Congr. image (mod 3)\\
    \hline 
     $\frac{23}{400}$   & 1200 &  NA  &   C  & $\{1\}$          &           &                  &  1\,852\,734\,273\,062\,400 \\
    \hline
  \end{tabular}
  \caption{Invariants for $\Gamma(5,1/5)$}
  \label{tab:m5-1.5}
\end{table}

\begin{table}[htbp]
  \centering
  {\large $\Gamma(6,1/6)$}\\
  \begin{tabular}[htbp]{|c|c|c|c|c|c|c|c|}
    \hline
   $\chi$     & LCM  & A? & C? & Ab               &  Cusps    & Cusp relations   & Congr. image (mod 5)\\
    \hline 
    $\frac{11}{144}$    &  144 &   NA &  NC  & $\Z_6$         & $R_1,R_2$ & $\br_3(R_1,R_2)$ & 152\,334\,000\,000  \\
    \hline
  \end{tabular}
  \begin{tabular}[htbp]{|c|c|c|c|c|c|}
    \hline
    Index                    & Core Index & Quotient                & Ab                                    & Self-int               & $b_1$ \\
    \hline
     $864=6\cdot 144$        &   20736    &                         & $\Z_3\oplus\Z_6^2\oplus\Z^2$          & $-1,-2,(-3)^5,(-6)^5$  & $2$\\
    \hline
  \end{tabular}
  \caption{Neat subgroup of $\Gamma(6,1/6)$ obtained from reduction mod 6}
  \label{tab:neat-m6-1.6}
\end{table}

\begin{table}[htbp]
  \centering
  {\large $\Gamma(8,1/8)$}\\
  \begin{tabular}[htbp]{|c|c|c|c|c|c|c|c|}
    \hline
   $\chi$     & LCM  & A? & C? & Ab               &  Cusps    & Cusp relations   & Congr. image (mod 3)\\
    \hline 
    $\frac{3}{32}$      &  192 &   A  &  C   & $\Z_8$         &           &                  & 42\,456\,960, $PGL(3,9)$\\
    \hline
  \end{tabular}
  \caption{Invariants for $\Gamma(8,1/8)$}
  \label{tab:m8-1.8}
\end{table}

\begin{table}[htbp]
  \centering
  {\large $\Gamma(12,1/12)$}\\
  \begin{tabular}[htbp]{|c|c|c|c|c|c|c|c|}
    \hline
   $\chi$     & LCM  & A? & C? & Ab               &  Cusps    & Cusp relations   & Congr. image (mod 5)\\
    \hline 
     $\frac{7}{72}$     & 144  &   A  &  C   & $\Z_{12}$        &           &                  & 152\,334\,000\,000  \\
    \hline
  \end{tabular}
  \caption{Invariants for $\Gamma(12,1/12)$}
  \label{tab:m12-1.12}
\end{table}

\begin{table}[htbp]
  \centering
  {\large $\Gamma(3,7/30)$}\\
  \begin{tabular}[htbp]{|c|c|c|c|c|c|c|c|}
    \hline
   $\chi$     & LCM  & A? & C? & Ab               &  Cusps    & Cusp relations   & Congr. image (mod 7)\\
    \hline 
     $\frac{2}{225}$    & 1800 &  A   &  C   & $\Z_3$         &           &                  &  ? \\
    \hline
  \end{tabular}
  \caption{Invariants for $\Gamma(3,7/30)$}
  \label{tab:m3-7.30}
\end{table}

\begin{table}[htbp]
  \centering
  {\large $\Gamma(4,3/20)$}\\
  \begin{tabular}[htbp]{|c|c|c|c|c|c|c|c|}
    \hline
   $\chi$     & LCM  & A? & C? & Ab               &  Cusps    & Cusp relations   & Congr. image (mod 3)\\
    \hline 
     $\frac{33}{800}$   & 2400 &  NA  &   C  & $\Z_2$         &           &                  & 1\,852\,734\,273\,062\,400  \\
    \hline
  \end{tabular}
  \caption{Invariants for $\Gamma(4,3/20)$}
  \label{tab:m4-3.20}
\end{table}

\begin{table}[htbp]
  \centering
  {\large $\Gamma(5,1/10)$}\\
  \begin{tabular}[htbp]{|c|c|c|c|c|c|c|c|}
    \hline
   $\chi$     & LCM  & A? & C? & Ab               &  Cusps    & Cusp relations   & Congr. image (mod 3)\\
    \hline 
    $\frac{13}{200}$    & 600  &   A  &  C   & $\Z_5$         &           &                  & 42\,573\,600, $PGU(3,9)$  \\
    \hline
  \end{tabular}
  \begin{tabular}[htbp]{|c|c|c|c|c|c|}
    \hline
    Index                    & Core Index & Quotient                & Ab                                    & Self-int            & $b_1$ \\
    \hline
     $600$                   & 600        &                         & $\Z^8$                                &                     & $8$\\
    \hline
     $600$                   & 15\,000     &                         & $\Z_5^2\oplus\Z^4$                    &                     & $4$\\
    \hline
     $600$                   & 15\,000     &                         & $\Z_5^5$                              &                     & $0$\\
    \hline
  \end{tabular}
 \caption{Neat subgroups of $\Gamma(5,1/10)$ of core index $\leq 30,000$}
 \label{tab:neat-m5-1.10}
\end{table}

\begin{table}[htbp]
  \centering
  {\large $\Gamma(10,0)$}\\
  \begin{tabular}[htbp]{|c|c|c|c|c|c|c|c|}
    \hline
   $\chi$     & LCM  & A? & C? & Ab               &  Cusps    & Cusp relations   & Congr. image (mod 3)\\
    \hline 
    $\frac{1}{10}$      & 150  &  A   & C    & $\Z_{10}$        &           &                  &  42\,573\,600, $PGU(3,9)$ \\
    \hline
  \end{tabular}
  \begin{tabular}[htbp]{|c|c|c|c|c|c|}
    \hline
    Index                    & Core Index & Quotient                & Ab                                    & Self-int            & $b_1$ \\
    \hline
     $150$                   &    600     &                         & $\Z_2^4\oplus\Z^4$                    &                     & $4$\\
    \hline
     $150$                   &    750     &                         & $\Z_5^2\oplus\Z^4$                    &                     & $4$\\
    \hline
     $150$                   &    1800    &                         & $\Z_5^2\oplus\Z^4$                    &                     & $4$\\
    \hline
  \end{tabular}
  \caption{Neat subgroups of $\Gamma(10,0)$ of core index $\leq 2,000$}
  \label{tab:neat-m10-0}
\end{table}

\begin{table}[htbp]
  \centering
  {\large $\Gamma(3,1/6)$}\\
  \begin{tabular}[htbp]{|c|c|c|c|c|c|c|c|}
    \hline
   $\chi$     & LCM  & A? & C? & Ab               &  Cusps    & Cusp relations   & Congr. image (mod 5)\\
    \hline 
     $\frac{1}{72}$     & 72   &  A   &  NC  & $\Z_3^2$         & {\scriptsize $c_1=1,c_2=J\bar2\bar1$} & {\scriptsize $c_1^3, c_2^6, \br_4(c_1,c_2)$} & 378\,000, $PGU(3,5)$  \\
    \hline
  \end{tabular}
  \begin{tabular}[htbp]{|c|c|c|c|c|c|}
    \hline
    Index                    & Core Index & Quotient                & Ab                                    & Self-int            & $b_1$ \\
    \hline
     $72$                    &    72      &                         & $\Z^4$                                &   $(-1)^4$          & $4$\\
    \hline
     $72$                    &   648      &                         & $\Z_3\oplus\Z^2$                      &   $-1,-3$           & $2$\\
    \hline
  \end{tabular}
  \caption{Neat subgroups of $\Gamma(3,1/6)$ of core index $\leq 2,000$}
  \label{tab:neat-m3-1.6}
\end{table}

\begin{table}[htbp]
  \centering
  {\large $\Gamma(4,1/12)$}\\
  \begin{tabular}[htbp]{|c|c|c|c|c|c|c|c|}
    \hline
   $\chi$     & LCM  & A? & C? & Ab               &  Cusps    & Cusp relations   & Congr. image (mod 5)\\
    \hline 
     $\frac{13}{288}$   & 288  &  NA  &  C   & $\Z_{12}$        &           &                  & 152\,334\,000\,000  \\
    \hline
  \end{tabular}
  \begin{tabular}[htbp]{|c|c|c|c|c|c|}
    \hline
    Index                    & Core Index & Quotient                & Ab                                    & Self-int            & $b_1$ \\
    \hline
     $864$                   &   864      &                         & $\Z_2^6\oplus\Z^4$                    &                     & $4$\\
    \hline
     $864$                   &  3456      &                         & $\Z_2^8\oplus\Z^2$                    &                     & $2$\\
    \hline
  \end{tabular}
  \caption{Neat subgroups of $\Gamma(4,1/12)$ of core index $\leq 10,000$}
  \label{tab:neat-m4-1.12}
\end{table}

\begin{table}[htbp]
  \centering
  {\large $\Gamma(6,0)$}\\
  \begin{tabular}[htbp]{|c|c|c|c|c|c|c|c|}
    \hline
   $\chi$     & LCM  & A? & C? & Ab               &  Cusps    & Cusp relations   & Congr. image (mod 5)\\
    \hline 
    $\frac{1}{12}$      &  36  &  A   &  NC  & $\Z_3\oplus\Z_6$ & $R_1,R_2$ & $\br_3(R_1,R_2)$ & 378\,000, $PGU(3,5)$  \\
    \hline
  \end{tabular}
  \begin{tabular}[htbp]{|c|c|c|c|c|c|}
    \hline
    Index                    & Core Index & Quotient                & Ab                                    & Self-int            & $b_1$ \\
    \hline
     $72$                    &    216     &                         & $\Z_3\oplus\Z^4$                      &  $(-1)^3,(-3)^3$    & $4$\\
    \hline
     $72$                    &    216     &                         & $\Z_3^2\oplus\Z^2$                    &  $(-3)^4$           & $2$\\
    \hline
     $72$                    &    648     &                         & $\Z_3\oplus\Z^2$                      &  $-3,-9$            & $2$\\
    \hline
     $72$                    &    864     &                         & $\Z_3\oplus\Z^2$                      &  $-12$              & $2$\\
    \hline
  \end{tabular}
  \caption{Neat subgroups of $\Gamma(6,0)$ of core index $\leq 1,000$}
  \label{tab:neat-m6-0}
\end{table}

\begin{table}[htbp]
  \centering
  {\large $\Gamma(3,5/42)$}\\
  \begin{tabular}[htbp]{|c|c|c|c|c|c|c|c|}
    \hline
   $\chi$     & LCM  & A? & C? & Ab               &  Cusps    & Cusp relations   & Congr. image (mod 5)\\
    \hline 
     $\frac{61}{3528}$  & 3528 &  NA  &  C   & $\Z_3$         &           &                  &  ? \\
    \hline
  \end{tabular}
  \caption{Invariants for $\Gamma(3,5/42)$}
  \label{tab:m3-5.42}
\end{table}

\begin{table}[htbp]
  \centering
  {\large $\Gamma(3,1/12)$}\\
  \begin{tabular}[htbp]{|c|c|c|c|c|c|c|c|}
    \hline
   $\chi$     & LCM  & A? & C? & Ab               &  Cusps    & Cusp relations   & Congr. image (mod 5)\\
    \hline 
     $\frac{11}{576}$   &  576 &  NA  &  C   & $\Z_3$         &           &                  & ?  \\
    \hline
  \end{tabular}
  \caption{Invariants for $\Gamma(3,1/12)$}
  \label{tab:m3-1.12}
\end{table}

\begin{table}[htbp]
  \centering
  {\large $\Gamma(4,0)$}\\
  \begin{tabular}[htbp]{|c|c|c|c|c|c|c|c|}
    \hline
   $\chi$     & LCM  & A? & C? & Ab               &  Cusps    & Cusp relations   & Congr. image (mod 3)\\
    \hline 
    $\frac{3}{64}$      & 192  &  A   & C    & $\Z_4$         &           &                  &  42\,456\,960, $PGL(3,9)$ \\
    \hline
  \end{tabular}
  \begin{tabular}[htbp]{|c|c|c|c|c|c|}
    \hline
    Index                    & Core Index            & Quotient                & Ab                                    & Self-int            & $b_1$ \\
    \hline
     $16512=86\cdot 192$     &  5663616 ($PGU(3,7)$) &                         & $?$                                    &                     & $?$\\
    \hline
  \end{tabular}
  \caption{Neat subgroups of $\Gamma(4,0)$}
  \label{tab:neat-m4-0}
\end{table}

\begin{table}[htbp]
  \centering
  {\large $\Gamma(3,1/18)$}\\
  \begin{tabular}[htbp]{|c|c|c|c|c|c|c|c|}
    \hline
   $\chi$     & LCM  & A? & C? & Ab               &  Cusps    & Cusp relations   & Congr. image (mod 5)\\
    \hline 
    $\frac{13}{648}$   & 648  &  A   &  C   & $\Z_3^2$         &           &                  &  ? \\
    \hline
  \end{tabular}
  \caption{Invariants for $\Gamma(3,1/18)$}
  \label{tab:m3-1.18}
\end{table}

\begin{table}[htbp]
  \centering
  {\large $\Gamma(3,1/30)$}\\
  \begin{tabular}[htbp]{|c|c|c|c|c|c|c|c|}
    \hline
   $\chi$     &  LCM  & A? & C? & Ab               &  Cusps    & Cusp relations   & Congr. image (mod 7)\\
    \hline 
    $\frac{37}{1800}$   & 1,800 &  NA  &  C   & $\Z_3$         &           &                  &  ? \\
    \hline
  \end{tabular}
  \caption{Invariants for $\Gamma(3,1/30)$}
  \label{tab:m3-1.30}
\end{table}

\begin{table}[htbp]
  \centering
  {\large $\Gamma(3,0)$}\\
  \begin{tabular}[htbp]{|c|c|c|c|c|c|c|c|}
    \hline
   $\chi$     & LCM  & A? & C? & Ab               &  Cusps    & Cusp relations   & Congr. image (mod 5)\\
    \hline 
     $\frac{1}{48}$     & 144  &  A   & C    & $\Z_3^2$         &           &                  &  ? \\
    \hline
  \end{tabular}
  \begin{tabular}[htbp]{|c|c|c|c|c|c|}
    \hline
    Index                    & Core Index & Quotient                & Ab                                    & Self-int            & $b_1$ \\
    \hline
     $2016=14\cdot 144$      &    54,432  &  $PGU(3,3)\times\Z_3^2$  & $\Z_3^4\oplus\Z^4$                    &                     & $4$\\
    \hline
     $2016$                  &    54,432  &  $PGU(3,3)\times\Z_3^2$  & $\Z_3^6\oplus\Z^2$                    &                     & $2$\\
    \hline
     $2016$                  &    54,432  &  $PGU(3,3)\times\Z_3^2$  & $\Z_3^8$                              &                     & $0$\\
    \hline
  \end{tabular}
  \caption{Neat subgroups of $\Gamma(3,0)$}
  \label{tab:neat-m3-0}
\end{table}

\end{document}